\documentclass[11pt]{amsart}
\usepackage[utf8]{inputenc}
\usepackage[T1]{fontenc}
\usepackage{amsmath,amssymb,amsthm,mathtools}
\usepackage[a4paper,margin=2.5cm]{geometry}
\usepackage{enumitem}

\usepackage{xcolor}
\usepackage[colorlinks=true,linkcolor=blue!60!black,citecolor=blue!60!black,urlcolor=blue!60!black]{hyperref}
\usepackage{orcidlink}

\makeatletter
\newcommand{\taggeditem}[2]{\item[#1]\phantomsection\protected@edef\@currentlabel{#1}\label{#2}}
\makeatother

\numberwithin{equation}{section}

\theoremstyle{plain}
\newtheorem{theorem}{Theorem}[section]
\newtheorem{proposition}[theorem]{Proposition}
\newtheorem{lemma}[theorem]{Lemma}
\theoremstyle{definition}
\newtheorem{definition}[theorem]{Definition}
\theoremstyle{remark}
\newtheorem{remark}[theorem]{Remark}

\newcommand{\R}{\mathbb{R}}
\newcommand{\Q}{\mathbb{Q}}
\newcommand{\N}{\mathbb{N}}
\newcommand{\dd}{\mathsf{d}}
\newcommand{\mm}{\mathfrak{m}}
\newcommand{\Ch}{\mathsf{Ch}}
\newcommand{\RCD}{\mathsf{RCD}}
\newcommand{\LIP}{\mathrm{LIP}}
\newcommand{\LIPbs}{\mathrm{LIP}_{\mathrm{bs}}}
\newcommand{\LIPb}{\mathrm{LIP}_{\mathrm{b}}}
\newcommand{\LIPfm}{\mathrm{LIP}_{\mathrm{fm}}}
\newcommand{\diff}{\mathrm{d}}
\newcommand{\PP}{\mathsf{P}}
\newcommand{\BE}{\mathsf{BE}}
\newcommand{\SL}{\mathsf{SL}}
\newcommand{\IH}{\mathsf{IH}}
\newcommand{\GVG}{\mathsf{G}}
\newcommand{\V}{\mathbb{V}}
\newcommand{\Test}{\mathrm{Test}^\infty(X)}
\newcommand{\Testp}{\mathrm{Test}^\infty_+(X)}
\newcommand{\G}{\Gamma}
\newcommand{\GG}{\mathbf{\Gamma}_2}

\DeclarePairedDelimiter{\abs}{\lvert}{\rvert}
\DeclarePairedDelimiter{\norm}{\lVert}{\rVert}

\title[Volume growth in the Eulerian definition of $\RCD$ spaces]{On the volume growth assumption in the Eulerian definition of $\RCD$ spaces}
\author[M. Picerni]{Marco Picerni \orcidlink{0009-0004-4364-4831}}
\address[Marco Picerni]{SISSA, Via Bonomea 265, 34136 Trieste, Italy}
\email{mpicerni@sissa.it}
\date{\today}
\subjclass[2020]{Primary 53C23; Secondary 31C25, 30L99, 35K05}
\keywords{RCD spaces, Bakry--\'Emery condition, volume growth, Radon measures, Dirichlet forms, Harnack inequality}

\begin{document}

\begin{abstract}
We prove that every infinitesimally Hilbertian metric measure space whose reference measure is a Radon measure, and which satisfies the Sobolev-to-Lipschitz property and the weak Bochner inequality with a lower Ricci bound, has Gaussian volume growth; in particular, balls have finite measure. Hence, in the Eulerian definition of $\RCD(K,\infty)$ spaces, both the volume growth assumption and the finiteness of the measure of balls are redundant. The proof combines the self-improvement of the Bakry--\'Emery gradient estimate, obtained without using the conservation of mass of the heat flow, with a dimension-free Harnack inequality derived from the Sobolev-to-Lipschitz property.
\end{abstract}

\maketitle

\section{Introduction}\label{sec:intro}

Lower bounds on the Ricci curvature of metric measure spaces were introduced, in a synthetic form based on optimal transport, by Sturm \cite{Stu06,Stu06b} and Lott--Villani \cite{LV09}, through the curvature-dimension condition $\mathsf{CD}(K,N)$; see also \cite{BS10,CM21} for the reduced condition $\mathsf{CD}^*(K,N)$ and its globalisation. The first-order calculus underlying the theory goes back to Cheeger \cite{Che99} and was developed, in the form used here, by Ambrosio, Gigli and Savar\'e \cite{AGS13,AGS14inv}, who also identified the heat flow as the $L^2$-gradient flow of the Cheeger energy and as the Wasserstein gradient flow of the entropy (see also \cite{Gig10}). Requiring in addition that the Cheeger energy be a quadratic form, that is, that the space be infinitesimally Hilbertian \cite{Gig15}, singles out the class of $\RCD(K,\infty)$ spaces \cite{AGS14b}, extended to $\sigma$-finite reference measures in \cite{AGMR15}; the finite-dimensional class $\RCD(K,N)$ was introduced in \cite{Gig15} and studied in \cite{EKS15,AMS19,CM21}. The Eulerian characterisation of these classes, via the weak Bochner inequality and the Sobolev-to-Lipschitz property, was obtained in \cite{AGS15} for $N=\infty$ and in \cite{EKS15,AMS19} for $N<\infty$; the classes are stable under (pointed) measured Gromov--Hausdorff convergence \cite{AGS14b,GMS15}. For the second-order calculus on these spaces we refer to \cite{Gig18} and to the monograph \cite{GP20}, whose notation we follow, and to the survey \cite{Amb18} for a general account of the theory.

In the Eulerian approach, $\RCD(K,\infty)$ spaces can be defined defined as follows (along the lines of \cite{AGS15}). Here $W^{1,2}(X)$, $\abs{Df}$, $\nabla f\cdot\nabla g$, $\Delta$ and $(h_t)_{t\geq0}$ denote the Sobolev space, the minimal weak upper gradient, the pointwise scalar product on the tangent module, the Laplacian and the heat flow; we recall them in Section~\ref{sec:prelim}.

\begin{definition}[$\RCD(K,\infty)$ spaces]\label{def:RCD}
Let $(X,\dd)$ be a complete and separable metric space, let $\mm$ be a Borel measure on $X$ which is finite on bounded sets and satisfies $\mathrm{supp}(\mm)=X$, and let $K\in\R$. The space $(X,\dd,\mm)$ is an $\RCD(K,\infty)$ space if the following hold:
\begin{itemize}[leftmargin=3.4em,labelsep=0.6em]
    \taggeditem{$(\GVG)$}{item:G} \emph{Gaussian volume growth}: there exist $\bar x\in X$ and $C,c>0$ such that $\mm(B_r(\bar x))\leq Ce^{cr^2}$ for every $r>0$;
    \taggeditem{$(\SL)$}{item:SL} \emph{Sobolev-to-Lipschitz property} \cite{AGS14b,Gig26}: every $f \in W^{1,2}(X)$ with $\abs{Df} \in L^\infty(\mm)$ admits a Lipschitz representative $\tilde{f} \in \LIP(X)$ such that $\tilde{f} = f$ $\mm$-a.e.\ and $\mathrm{Lip}(\tilde{f}) \leq \norm{\abs{Df}}_{L^\infty(\mm)}$;
    \taggeditem{$(\IH)$}{item:IH} \emph{Infinitesimal Hilbertianity} \cite{Gig15}: the Sobolev space $W^{1,2}(X)$ endowed with the norm $\norm{f}_{W^{1,2}(X)}^2 \coloneqq \norm{f}_{L^2(\mm)}^2 + 2\Ch(f)$, where $\Ch(f) \coloneqq \frac{1}{2}\int_X \abs{Df}^2 \diff\mm$, is a Hilbert space;
    \taggeditem{$(\BE)$}{item:BE} \emph{Weak Bochner inequality} $\BE(K,\infty)$: for every $f \in D(\Delta)$ with $\Delta f \in W^{1,2}(X)$ and every $\varphi \in D(\Delta) \cap L^\infty(\mm)$ with $\varphi \geq 0$ and $\Delta \varphi \in L^\infty(\mm)$,
    \begin{equation}\label{eq:bochner_weak}
        \frac{1}{2} \int_X \abs{Df}^2 \Delta \varphi \diff\mm
        \geq \int_X \varphi \bigl( \nabla f \cdot \nabla \Delta f + K \abs{Df}^2 \bigr) \diff\mm.
    \end{equation}
\end{itemize}
\end{definition}

\begin{remark}\label{rem:lagrangian}
In the Lagrangian approach, based on the $K$-convexity of the Boltzmann--Shannon entropy along Wasserstein geodesics, \ref{item:G} is a consequence of the curvature condition: Sturm proved that every $\mathsf{CD}(K,\infty)$ space satisfies $\int_Xe^{-c\,\dd^2(x_0,\cdot)}\diff\mm<\infty$ for some $c>0$ \cite[Thm.~4.24]{Stu06}, which gives \ref{item:G} since $\mm(B_r(x_0))\,e^{-cr^2}\leq\int_{B_r(x_0)}e^{-c\,\dd^2(x_0,\cdot)}\diff\mm$ (see also \cite[Section~2]{AGMR15}). In the Eulerian approach, \ref{item:G} is part of the definition, and it is used in \cite{AGS15} to prove the equivalence with the Lagrangian approach.
\end{remark}

The purpose of this note is to show that, in Definition~\ref{def:RCD}, \ref{item:G} follows from \ref{item:SL}, \ref{item:IH} and \ref{item:BE}, and that the finiteness of $\mm$ on bounded sets can be weakened to the requirement that $\mm$ be a Radon measure.

\begin{definition}[Radon measures]\label{def:radon}
Let $(X,\dd)$ be a metric space. A Borel measure $\mm$ on $X$ is a \emph{Radon measure} if $\mm(K)<\infty$ for every compact set $K\subseteq X$, $\mm$ is outer regular and $\mm$ is inner regular on open sets, that is,
\begin{gather*}
    \mm(B)=\inf\bigl\{\mm(U):U\supseteq B,\ U\text{ open}\bigr\}\qquad\text{for every Borel set }B\subseteq X,\\
    \mm(U)=\sup\bigl\{\mm(K):K\subseteq U,\ K\text{ compact}\bigr\}\qquad\text{for every open set }U\subseteq X .
\end{gather*}
\end{definition}

\medskip
\noindent\textbf{(S)} $(X,\dd)$ is a complete and separable metric space, $\mm$ is a Radon measure on $X$ with $\mathrm{supp}(\mm)=X$, and \ref{item:SL}, \ref{item:IH} and \ref{item:BE} hold for some $K\in\R$.
\medskip

The objects in \ref{item:SL}, \ref{item:IH} and \ref{item:BE} make sense under \textbf{(S)}; see Section~\ref{sec:prelim}. The regularity in Definition~\ref{def:radon} cannot be dropped: in Remark~\ref{rem:radon_examples} we give a Borel measure on a separable Hilbert space, finite on compact sets and with full support, for which \ref{item:SL}, \ref{item:IH} and \ref{item:BE} hold trivially and every ball has infinite measure. For $K\in\R$ and $t>0$ we set, following the notation of \cite{AGS15}
\begin{equation}\label{eq:cK}
    I_K(t) \coloneqq \int_0^t e^{2Ks} \diff s = \frac{e^{2Kt} - 1}{2K} \quad (\text{with } c_0(t) \coloneqq t).
\end{equation}

\begin{theorem}[Gaussian volume growth]\label{thm:main}
Assume \textbf{(S)}. Let $x_0\in X$, let $r_0>0$ be such that $\mm(\overline{B}_{r_0}(x_0))<\infty$, let $t>0$, and set $f\coloneqq\frac12\bigl(1-\dd(\cdot,x_0)/r_0\bigr)^+$. Then the continuous representative of $h_tf$ satisfies $h_tf(x_0)>0$, and for every $r>0$
\begin{equation}\label{eq:main_growth}
    \mm\bigl(B_r(x_0)\bigr) \leq \frac{\norm{f}_{L^1(\mm)}}{\bigl(h_t f(x_0)\bigr)^2} \exp\Bigl( \frac{r^2}{2 I_K(t)} \Bigr).
\end{equation}
In particular, $\mm$ is finite on bounded sets, and $(X,\dd,\mm)$ is an $\RCD(K,\infty)$ space.
\end{theorem}

A radius $r_0$ as in Theorem~\ref{thm:main} exists for every $x_0\in X$ by Lemma~\ref{lem:locally_finite}.

\begin{remark}[Dependence of the constants]\label{rem:constants}
Choosing $t = 1$ in \eqref{eq:main_growth} yields $\mm(B_r(x_0)) \leq C e^{c r^2}$ with $c = (2I_K(1))^{-1}$ and $C = \norm{f}_{L^1(\mm)}(h_1 f(x_0))^{-2}$. A quantitative upper bound on $C$ requires a lower bound on $h_t f(x_0)$, which can be obtained from the short-time estimate $\norm{h_t f - f}_{L^2(\mm)}^2 \leq t \norm{Df}_{L^2(\mm)}^2$ combined with the Harnack inequality of Lemma~\ref{lem:power_harnack}; this brings in $\mm(B_{r_0/2}(x_0))$ and makes the choice of $t$ (hence $c$) depend on $K$, $r_0$ and the ratio $\mm(\overline{B}_{r_0}(x_0))/\mm(B_{r_0/2}(x_0))$.
\end{remark}

\subsection{Strategy of the proof and related results}\label{subsec:strategy}

The proof follows a route that is classical for smooth weighted manifolds; the point is to check that, in the present nonsmooth setting, none of its steps requires the growth of $\mm$, nor its finiteness on balls.

The starting point is the $L^2$ Bakry--\'Emery estimate $\abs{Dh_tf}^2\leq e^{-2Kt}h_t(\abs{Df}^2)$, which is proved in \cite[Cor.~2.3]{AGS15} starting from the weak Bochner inequality; its proof uses neither \ref{item:G} nor the finiteness of $\mm$ on balls (Lemma~\ref{lem:GP}).

The second ingredient is the $L^1$ Bakry--\'Emery estimate $\abs{Dh_tf}\leq e^{-Kt}h_t(\abs{Df})$ (Proposition~\ref{prop:L1BE}), due to Bakry for smooth diffusions and to Savar\'e \cite{Sav14} for abstract Dirichlet forms. We follow \cite{Sav14}, but without requiring that the heat flow preserves mass. Savar\'e uses this assumption to construct the iterated carr\'e du champ as a finite measure \cite[Lemma~2.6]{Sav14}; we represent instead the nonnegative functional $\GG^K[f]$ by a locally finite measure (Lemma~\ref{lem:repr}). Moreover, the quasi-regularity of the Cheeger energy, proved in \cite[Thm.~4.1]{Sav14} by exhausting $X$ with balls, is obtained by exhausting $X$ with superlevel sets of test functions (Lemma~\ref{lem:QR}).

The third ingredient is a dimension-free Harnack inequality (Lemma~\ref{lem:power_harnack}): for $0\leq f\le1/2$,
\[
    \bigl(h_tf(x)\bigr)^2\leq h_tf(y)\exp\Bigl(\frac{\dd(x,y)^2}{2I_K(t)}\Bigr)\qquad\text{for all } x,y\in X.
\]
Wang \cite{Wan97} proved on Riemannian manifolds the stronger inequality
\[
    \bigl(h_tf(x)\bigr)^2\leq h_t(f^2)(y)\exp\Bigl(\frac{\dd(x,y)^2}{2I_K(t)}\Bigr),
\]
and Li \cite{Li16} extended it to $\RCD(K,\infty)$ spaces defined via optimal transport. Li's proof uses the $W_2$-contractivity of the heat flow, which relies on the identification of the heat flow with the Wasserstein gradient flow of the entropy and is not available under \textbf{(S)}. This is why we prove the weaker inequality above: it only requires the $L^1$ Bakry--\'Emery estimate and the Sobolev-to-Lipschitz property, and it suffices to bound the volume of balls.

The $L^1$ Bakry--\'Emery estimate is essential here. Indeed, for every $\delta\in(0,1/4)$, the proof of Lemma~\ref{lem:power_harnack} shows that $\sqrt{-\log(h_tf+\delta)}$ is $(2\sqrt{I_K(t)})^{-1}$-Lipschitz; this bounds $h_tf(y)$ from below at every $y\in X$ in terms of $h_tf(x)$ and $\dd(x,y)$. The key step is the pointwise bound
\[
    \abs{Dh_tf}^2\le e^{-2Ks}(h_tf+\delta)\,h_s\Bigl(\frac{\abs{Dh_{t-s}f}^2}{h_{t-s}f+\delta}\Bigr),
\]
which follows from $\abs{Dh_tf}\le e^{-Ks}h_s(\abs{Dh_{t-s}f})$ and the Cauchy--Schwarz inequality for $h_s$. The $L^2$ Bakry--\'Emery estimate $\abs{Dh_tf}^2\le e^{-2Ks}h_s(\abs{Dh_{t-s}f}^2)$ only gives the same bound with the factor $h_tf+\delta$ replaced by the constant $\norm{h_{t-s}f+\delta}_{L^\infty(\mm)}\le 1/2+\delta$. Following the rest of the argument, this yields that $G(h_tf+\delta)$ is $I_K(t)^{-1/2}$-Lipschitz, where $G(u)\coloneqq\int_0^u\bigl(v\log(1/v)\bigr)^{-1/2}\diff v$ for $u\in[0,3/4]$. Since $G$ is bounded, such an estimate gives no information on $h_tf(y)$ once $\dd(x,y)\ge I_K(t)^{1/2}G(3/4)$, and therefore no lower bound at every point.

Finally, Theorem~\ref{thm:main} follows by integrating the Harnack inequality in $y$ over $B_r(x_0)$ and using the $L^1$-contractivity of $h_t$.

For smooth weighted manifolds $(M,g,e^{-V}\mathrm{vol}_g)$ whose Bakry--\'Emery Ricci tensor satisfies
\[
    \mathrm{Ric}_V\coloneqq\mathrm{Ric}_g+\mathrm{Hess}_gV\geq K,
\]
with no assumption on the potential $V$, the Gaussian growth $(e^{-V}\mathrm{vol}_g)(B_r(x_0))\leq Ce^{cr^2}$ is due to Morgan \cite{Mor05} and Wei--Wylie \cite[Thm.~3.1]{WW09}, and it is sharp on the Gaussian soliton. Their proof relies on comparison geometry along geodesics, which has no direct counterpart in the Eulerian nonsmooth setting; the argument above replaces it by the heat flow. A discrete counterpart, for graphs of bounded degree satisfying the dimension-free Bakry--\'Emery condition with $K=0$, has recently been obtained in \cite{Graph26}.

\section{Preliminaries}\label{sec:prelim}

\subsection{Radon measures}\label{subsec:radon}

We will use the following elementary lemma.

\begin{lemma}[Local finiteness]\label{lem:locally_finite}
Let $(X,\dd)$ be a metric space and let $\mm$ be a Radon measure on $X$. Then every $x\in X$ has a neighbourhood of finite measure. If moreover $(X,\dd)$ is separable, then there are $x_i\in X$ and $r_i>0$, $i\in\N$, such that $\mm(B_{3r_i}(x_i))<\infty$ for every $i$ and $X=\bigcup_iB_{r_i}(x_i)$; in particular, $\mm$ is $\sigma$-finite and $L^2(\mm)$ is separable.
\end{lemma}

\begin{proof}
Let $x\in X$. Since $\{x\}$ is compact, $\mm(\{x\})<\infty$, and by outer regularity there is an open set $U\ni x$ with $\mm(U)<\infty$; $U$ contains a ball centred at $x$. If $(X,\dd)$ is separable, it is Lindel\"of, so the open cover $\{B_{r_x}(x)\}_{x\in X}$, where $r_x>0$ is such that $\mm(B_{3r_x}(x))<\infty$, has a countable subcover. The separability of $L^2(\mm)$ follows from the $\sigma$-finiteness of $\mm$ and the separability of $(X,\dd)$.
\end{proof}

Under \textbf{(S)}, we call a ball $B_r(x)$ \emph{admissible} if $\mm(B_{3r}(x))<\infty$. By Lemma~\ref{lem:locally_finite}, every point is the centre of admissible balls, and $X$ is covered by countably many admissible balls.

\begin{remark}\label{rem:radon_examples}
On a complete and separable metric space, a Borel measure is a Radon measure if and only if every point has a neighbourhood of finite measure. One implication is Lemma~\ref{lem:locally_finite}. Conversely, if $X$ is covered by countably many open sets $U_i$ of finite measure, then $\mm$ is finite on compact sets, and outer regularity and inner regularity on open sets follow from the corresponding properties of the finite measures $\mm\llcorner U_i$, the latter by Ulam's theorem.

The finiteness on compact sets alone is not enough for Theorem~\ref{thm:main}. Let $X$ be a separable infinite-dimensional Hilbert space. Compact subsets of $X$ have empty interior, so $\sigma$-compact subsets of $X$ are meagre, and no nonempty open set is contained in a $\sigma$-compact set by the Baire category theorem. Set $\mm(B)\coloneqq0$ if the Borel set $B$ is contained in a $\sigma$-compact set and $\mm(B)\coloneqq\infty$ otherwise. Then $\mm$ is a Borel measure (if one of the disjoint sets $B_k$ is not contained in a $\sigma$-compact set, neither is $\bigcup_kB_k$), it vanishes on compact sets, and $\mm(U)=\infty$ for every nonempty open set $U$, so $\mathrm{supp}(\mm)=X$. Every $f\in L^2(\mm)$ vanishes $\mm$-a.e., since the sets $\{\abs{f}>1/k\}$ have finite, hence zero, measure. Thus $W^{1,2}(X)=\{0\}$, and \ref{item:SL}, \ref{item:IH} and \ref{item:BE} hold trivially, while every ball has infinite measure. This measure is neither outer regular, since $\mm(\{0\})=0$ while every open set containing $0$ has infinite measure, nor inner regular on open sets.
\end{remark}

We denote by $\LIPfm(X)$ the space of bounded Lipschitz functions $f\colon X\to\R$ such that $\mm(\mathrm{supp}(f))<\infty$. Under \textbf{(S)}, this space replaces the space $\LIPbs(X)$ of Lipschitz functions with bounded support, which in general is not contained in $L^2(\mm)$.

\begin{lemma}\label{lem:LIPfm}
Assume that $(X,\dd)$ is complete and separable and that $\mm$ is a Radon measure on $X$. Then:
\begin{enumerate}[label=(\roman*),ref=(\roman*),leftmargin=*]
    \item\label{item:LIPfm_sob} $\LIPfm(X)\subseteq W^{1,2}(X)$, and $\abs{Df}\leq\mathrm{lip}(f)\leq\mathrm{Lip}(f)\,\chi_{\mathrm{supp}(f)}$ $\mm$-a.e.\ for every $f\in\LIPfm(X)$;
    \item\label{item:LIPfm_bump} if $B_r(x)$ is admissible, then $\beta_{x,r}\coloneqq\bigl(1-\dd(\cdot,B_r(x))/r\bigr)^+$ belongs to $\LIPfm(X)$, $0\leq\beta_{x,r}\leq1$, $\beta_{x,r}=1$ on $B_r(x)$ and $\beta_{x,r}=0$ outside $B_{2r}(x)$;
    \item\label{item:LIPfm_dense} $\LIPfm(X)$ is dense in $L^2(\mm)$.
\end{enumerate}
\end{lemma}

\begin{proof}
\ref{item:LIPfm_sob}. Since $f$ is bounded and $\mm(\mathrm{supp}(f))<\infty$, $f\in L^2(\mm)$. The local Lipschitz constant $\mathrm{lip}(f)$ vanishes on the open set $X\setminus\mathrm{supp}(f)$ and is bounded by $\mathrm{Lip}(f)$, so $\mathrm{lip}(f)\in L^2(\mm)$. Moreover, $\mathrm{lip}(f)$ is an upper gradient of $f$: for every absolutely continuous curve $\gamma\colon[0,1]\to X$, the function $t\mapsto f(\gamma_t)$ is absolutely continuous and $\abs{f(\gamma_1)-f(\gamma_0)}\leq\int_0^1\mathrm{lip}(f)(\gamma_t)\,\abs{\dot\gamma_t}\diff t$. In particular, $\mathrm{lip}(f)$ is a weak upper gradient of $f$, so $f\in W^{1,2}(X)$ and $\abs{Df}\leq\mathrm{lip}(f)$ $\mm$-a.e.

\ref{item:LIPfm_bump}. $\beta_{x,r}$ is $r^{-1}$-Lipschitz, and its support is contained in $\overline{B}_{2r}(x)\subseteq B_{3r}(x)$.

\ref{item:LIPfm_dense}. Let $(B_{r_i}(x_i))_{i\in\N}$ be admissible balls covering $X$, and set $E_n\coloneqq\bigcup_{i\leq n}B_{r_i}(x_i)$. Then $E_n$ is open, $E_n\uparrow X$, and $\mm(\overline{E}_n)\leq\sum_{i\leq n}\mm(B_{3r_i}(x_i))<\infty$. Let $g\in L^2(\mm)$ and $\varepsilon>0$. By dominated convergence there is $n$ with $\norm{g-g\chi_{E_n}}_{L^2(\mm)}<\varepsilon$. Since bounded Lipschitz functions are dense in $L^2(\nu)$ for every finite Borel measure $\nu$ on a metric space, applying this to $\nu\coloneqq\mm\llcorner E_n$ we find a bounded Lipschitz function $u$ with $\norm{u-g}_{L^2(\mm\llcorner E_n)}<\varepsilon$. The functions $\zeta_k\coloneqq\min\{1,k\,\dd(\cdot,X\setminus E_n)\}$ are Lipschitz, vanish outside $E_n$ and satisfy $\zeta_k\uparrow\chi_{E_n}$, so $\zeta_ku\in\LIPfm(X)$, because $\mathrm{supp}(\zeta_ku)\subseteq\overline{E}_n$, and $\norm{\zeta_ku-u}_{L^2(\mm\llcorner E_n)}\to0$ by dominated convergence. Hence $\norm{\zeta_ku-g}_{L^2(\mm)}<3\varepsilon$ for $k$ large.
\end{proof}

\subsection{Calculus and heat flow}\label{subsec:calculus}

Under \textbf{(S)}, the Sobolev space $W^{1,2}(X)$ and the minimal weak upper gradient $\abs{Df}$ are defined via test plans; under \ref{item:IH}, $(f,g)\mapsto\nabla f\cdot\nabla g\in L^1(\mm)$ is the pointwise scalar product on the Hilbert tangent module $L^2(TX)$, $\Delta\colon D(\Delta)\subseteq W^{1,2}(X)\to L^2(\mm)$ is the Laplacian associated with $\Ch$, and $(h_t)_{t\geq0}$ is the heat flow, i.e.\ the gradient flow of $\Ch$ in $L^2(\mm)$.

\begin{proposition}\label{prop:GP_calculus}
Let $(X,\dd)$ be a complete and separable metric space, let $\mm$ be a Radon measure on $X$ with $\mathrm{supp}(\mm)=X$, and assume \ref{item:IH}. Then:
\begin{enumerate}[label=(\alph*),ref=(\alph*),leftmargin=*]
    \item\label{item:calc_rules} \emph{Calculus rules}: the minimal weak upper gradient is local; for $f\in W^{1,2}(X)$ and $\eta\colon\R\to\R$ Lipschitz with $\eta(0)=0$, $\eta\circ f\in W^{1,2}(X)$ and $\abs{D(\eta\circ f)}=\abs{\eta'}\circ f\,\abs{Df}$ $\mm$-a.e.; $\nabla f\cdot\nabla g$ is bilinear, symmetric and continuous from $W^{1,2}(X)\times W^{1,2}(X)$ to $L^1(\mm)$, and satisfies the chain rule and the Leibniz rule on $W^{1,2}(X)\cap L^\infty(\mm)$; if $f_n\to f$ in $L^2(\mm)$, $f_n\to f$ $\mm$-a.e.\ and $\abs{Df_n}\rightharpoonup G$ weakly in $L^2(\mm)$, then $G$ is a weak upper gradient of $f$.
    \item\label{item:calc_lap} \emph{Laplacian}: $\Delta$ is self-adjoint, and
    \[
        \int_Xg\,\Delta f\diff\mm=-\int_X\nabla f\cdot\nabla g\diff\mm\qquad\text{for } f\in D(\Delta),\ g\in W^{1,2}(X);
    \]
    for $f,g\in\LIPb(X)\cap D(\Delta)$ and $\phi\in C^2(\R)$ with $\phi''$ bounded,
    \[
        \Delta(fg)=f\Delta g+g\Delta f+2\nabla f\cdot\nabla g,\qquad\Delta(\phi\circ f)=\phi'\circ f\,\Delta f+\phi''\circ f\,\abs{Df}^2 .
    \]
    \item\label{item:calc_heat} \emph{Heat flow}: $h_t$ is linear and self-adjoint, $\Ch(h_tf)\leq\Ch(f)$, $h_t\Delta f=\Delta h_tf$ for $f\in D(\Delta)$; for $f\in L^2(\mm)$, $t\mapsto h_tf$ belongs to $C([0,\infty);L^2(\mm))\cap C^\infty((0,\infty);W^{1,2}(X))$, with $h_tf\in D(\Delta)$ and $\frac{\diff}{\diff t}h_tf=\Delta h_tf$ for $t>0$; $(h_t)_{t\geq0}$ is positivity preserving, satisfies the weak maximum principle, and extends to a semigroup of contractions on $L^p(\mm)$ for every $p\in[1,\infty]$.
\end{enumerate}
\end{proposition}

\begin{proof}
Originally proved in \cite{AGS14b,Gig15}, all these statements are also presented in \cite[Chapters~2, 4, and~5]{GP20} under the standing assumption (1.21) therein that $\mm$ is finite on bounded sets; we shall follow the latter presentation throughout. In Chapters~1--5 of \cite{GP20}, this assumption is used explicitly only in Prop.~1.1.18 (density of $\LIPbs(X)$ in $L^p$), Prop.~1.1.19 (the topology of $L^0(\mm)$, where $\sigma$-finiteness suffices), Prop.~2.1.26 ($\LIPbs(X)\subseteq W^{1,2}(X)$) and Thm.~2.2.28 (density in energy of Lipschitz functions). In our setting, $\mm$ is $\sigma$-finite by Lemma~\ref{lem:locally_finite}, and wherever the inclusion $\LIPbs(X)\subseteq W^{1,2}(X)$ or the density of $\LIPbs(X)$ in $L^2(\mm)$ is used, Lemma~\ref{lem:LIPfm} provides the same statements for $\LIPfm(X)$; in particular $W^{1,2}(X)$ is dense in $L^2(\mm)$, which is needed for the heat flow to be defined on the whole $L^2(\mm)$. Thm.~2.2.28 is not used in this note.
\end{proof}

\subsection{Growth-free results on the heat flow}\label{subsec:GP61}

The results of this subsection are proved in \cite{AGS15,Sav14}. Following  the proofs collected in \cite[Chapters~1--6]{GP20} for $\RCD(K,\infty)$ spaces, one can check that they only use the results of Chapters~1--5, which hold in our setting by Proposition~\ref{prop:GP_calculus}, together with \ref{item:SL}, \ref{item:IH} and \ref{item:BE}. In particular, they use neither \ref{item:G} nor the finiteness of $\mm$ on balls of arbitrary radius.

We use the space of test functions
\[
    \Test \coloneqq \bigl\{ f \in \LIP(X) \cap L^\infty(\mm) \cap D(\Delta) : \Delta f \in W^{1,2}(X) \cap L^\infty(\mm) \bigr\},
\]
its subset $\Testp$ of $\mm$-a.e.\ nonnegative functions, and, for $f \in L^2(\mm)$ and $\varphi \in C_c^\infty(0,\infty)$, the mollified heat flow
\[
    h_\varphi f \coloneqq \int_0^\infty h_s f\,\varphi(s)\diff s .
\]
Elements of $\Test$ are identified with their Lipschitz representatives, which are unique since $\mathrm{supp}(\mm)=X$.

\begin{lemma}[Growth-free results]\label{lem:GP}
Assume \textbf{(S)}. Then:
\begin{enumerate}[label=(\roman*),ref=(\roman*),leftmargin=*]
    \item\label{item:GP_BE} for every $f \in W^{1,2}(X)$ and every $t \geq 0$,
    \begin{equation}\label{eq:BE_quadratic}
        \abs{D h_t f}^2 \leq e^{-2Kt} h_t\bigl(\abs{Df}^2\bigr) \quad \mm\text{-a.e.\ in } X;
    \end{equation}
    \item\label{item:GP_Lip} for every $f \in L^2(\mm) \cap L^\infty(\mm)$ and every $t > 0$, $\abs{Dh_tf} \in L^\infty(\mm)$ and $h_t f$ admits a Lipschitz representative;
    \item\label{item:GP_test} for every $f \in L^2(\mm) \cap L^\infty(\mm)$ and every $\varphi \in C_c^\infty(0,\infty)$, $h_\varphi f \in \Test$;
    \item\label{item:GP_dense} $\Testp$ is dense in $W^{1,2}(X)_+$ and $\Test$ is dense in $W^{1,2}(X)$;
    \item\label{item:GP_sq} $\abs{Df}^2\in W^{1,2}(X)$ for every $f\in\Test$;
    \item\label{item:GP_alg} $\Test$ is an algebra.
\end{enumerate}
\end{lemma}

\section{Dirichlet forms}\label{sec:dirichlet}

We collect the notions from the theory of Dirichlet forms used in Section~\ref{sec:proof}, following \cite[Sections~2--3]{Sav14} and \cite[Section~2]{AGS15}. Throughout this section, $(X,\tau)$ is a Polish space, $\mm$ is a $\sigma$-finite Borel measure on $X$ with $\mathrm{supp}(\mm)=X$, and functions are identified when they agree $\mm$-a.e.

\begin{definition}[Symmetric Dirichlet form, generator, semigroup]\label{def:dirichlet}
A \emph{symmetric Dirichlet form} on $L^2(X,\mm)$ is a lower semicontinuous quadratic form $\mathcal{E}\colon L^2(X,\mm)\to[0,\infty]$ whose domain $\V\coloneqq\{f\in L^2(X,\mm):\mathcal{E}(f)<\infty\}$ is dense in $L^2(X,\mm)$ and which satisfies the Markov property
\[
    \mathcal{E}(\eta\circ f)\le\mathcal{E}(f)\qquad\text{for every } f\in L^2(X,\mm) \text{ and every $1$-Lipschitz } \eta\colon\R\to\R \text{ with } \eta(0)=0 .
\]
We denote by $\mathcal{E}(f,g)\coloneqq\frac14\bigl(\mathcal{E}(f+g)-\mathcal{E}(f-g)\bigr)$ the associated symmetric bilinear form on $\V$, by $\norm{f}_{\V}^2\coloneqq\norm{f}_{L^2(\mm)}^2+\mathcal{E}(f)$ the Hilbert norm of $\V$, and set $\V_\infty\coloneqq\V\cap L^\infty(\mm)$ and $\V_+\coloneqq\{f\in\V:f\ge0\}$. The \emph{generator} of $\mathcal{E}$ is the nonpositive self-adjoint operator $\mathsf{L}\colon D(\mathsf{L})\subseteq\V\to L^2(X,\mm)$ characterised by
\[
    f\in D(\mathsf{L}),\ \mathsf{L}f=h \iff f\in\V \text{ and } \mathcal{E}(f,g)=-\int_X hg\diff\mm \text{ for every } g\in\V,
\]
and $(\PP_t)_{t\ge0}$, $\PP_t\coloneqq e^{t\mathsf{L}}$, is the associated \emph{Markov semigroup}.
\end{definition}

\begin{remark}\label{rem:semigroup}
By the Markov property, $(\PP_t)_{t\ge0}$ is positivity preserving and sub-Markovian, that is, $0\leq f\leq c$ $\mm$-a.e.\ implies $0\le\PP_tf\leq c$ $\mm$-a.e.; it extends to a semigroup of contractions on $L^p(X,\mm)$ for every $p\in[1,\infty]$, strongly continuous for $p<\infty$ and weakly$^*$ continuous for $p=\infty$, and it is an analytic semigroup on $L^2(X,\mm)$ \cite[Section~2.1]{AGS15}.
\end{remark}

\begin{definition}[Strong locality and carr\'e du champ]\label{def:local_cdc}
Let $\mathcal{E}$ be a symmetric Dirichlet form on $L^2(X,\mm)$.
\begin{enumerate}[label=(\roman*),ref=(\roman*),leftmargin=*]
    \item $\mathcal{E}$ is \emph{strongly local} if $\mathcal{E}(f,g)=0$ for all $f,g\in\V$ such that $(f+a)g=0$ $\mm$-a.e.\ for some $a\in\R$.
    \item $\mathcal{E}$ \emph{admits a carr\'e du champ} if there is a symmetric, bilinear and continuous map $\Gamma\colon\V\times\V\to L^1(X,\mm)$ such that
    \[
        2\int_X\Gamma(f,g)\,\varphi\diff\mm=\mathcal{E}(f,g\varphi)+\mathcal{E}(g,f\varphi)-\mathcal{E}(fg,\varphi)\qquad\text{for all } f,g,\varphi\in\V_\infty .
    \]
    We write $\Gamma(f)\coloneqq\Gamma(f,f)$.
\end{enumerate}
\end{definition}

\begin{definition}[Capacity, nests, quasi-continuity]\label{def:capacity}
Let $\mathcal{E}$ be a symmetric Dirichlet form on $L^2(X,\mm)$. The \emph{capacity} of an open set $A\subseteq X$ is
\[
    \mathrm{Cap}(A)\coloneqq\inf\bigl\{\norm{u}_{\V}^2:u\in\V,\ u\ge1\ \mm\text{-a.e.\ on } A\bigr\}\qquad(\inf\emptyset\coloneqq+\infty),
\]
and the capacity of an arbitrary set $B\subseteq X$ is $\mathrm{Cap}(B)\coloneqq\inf\{\mathrm{Cap}(A):A\supseteq B,\ A \text{ open}\}$. For a closed set $F\subseteq X$ let $\V_F\coloneqq\{f\in\V:f=0\ \mm\text{-a.e.\ on } X\setminus F\}$.
\begin{enumerate}[label=(\roman*),ref=(\roman*),leftmargin=*]
    \item An \emph{$\mathcal{E}$-nest} is a nondecreasing sequence $(F_k)_{k\in\N}$ of closed subsets of $X$ such that $\bigcup_k\V_{F_k}$ is dense in $\V$.
    \item A set $N\subseteq X$ is \emph{$\mathcal{E}$-polar} if $N\subseteq X\setminus\bigcup_kF_k$ for some $\mathcal{E}$-nest $(F_k)_{k\in\N}$; a property holds \emph{quasi-everywhere} (q.e.) if it holds outside an $\mathcal{E}$-polar set.
    \item A function $f\colon X\to\R$ is \emph{$\mathcal{E}$-quasi-continuous} if there is an $\mathcal{E}$-nest $(F_k)_{k\in\N}$ such that the restriction $f|_{F_k}$ is continuous on $F_k$, endowed with the topology induced by $X$, for every $k$.
\end{enumerate}
\end{definition}

\begin{lemma}\label{lem:sep_cap}
Let $\mathcal{E}$ be a symmetric Dirichlet form on $L^2(X,\mm)$. Then $\V$ is separable, and $\mathrm{Cap}\bigl(\bigcup_kA_k\bigr)\leq\sum_k\mathrm{Cap}(A_k)$ for every sequence $(A_k)_{k\in\N}$ of subsets of $X$.
\end{lemma}

\begin{proof}
The countable subadditivity of the capacity holds for every symmetric Dirichlet form \cite[Ch.~III, Thm.~2.8]{MR92}. The space $\V$, endowed with $\norm{\cdot}_\V$, is a Hilbert space, since $\mathcal{E}$ is lower semicontinuous, and the inclusion $\iota\colon\V\to L^2(\mm)$ is linear, continuous and injective; hence its adjoint $\iota^*\colon L^2(\mm)\to\V$ has dense range, since $(\mathrm{ran}\,\iota^*)^\perp=\ker\iota=\{0\}$. As $X$ is Polish and $\mm$ is $\sigma$-finite, $L^2(\mm)$ is separable, hence so are $\mathrm{ran}\,\iota^*$ and $\V$.
\end{proof}

\begin{definition}[Quasi-regularity]\label{def:QR}
A symmetric Dirichlet form $\mathcal{E}$ on $L^2(X,\mm)$ is \emph{quasi-regular} if:
\begin{enumerate}[label=$\langle\mathrm{QR}.\arabic*\rangle$,ref=$\langle\mathrm{QR}.\arabic*\rangle$,leftmargin=*]
    \item there exists an $\mathcal{E}$-nest consisting of compact sets;
    \item there exists a dense subset of $\V$ whose elements admit $\mathcal{E}$-quasi-continuous representatives;
    \item there exist an $\mathcal{E}$-polar set $N\subseteq X$ and a countable family of $\mathcal{E}$-quasi-continuous functions in $\V$ separating the points of $X\setminus N$.
\end{enumerate}
\end{definition}

\begin{lemma}[{\cite[Lemma~2.4]{Sav14}}]\label{lem:Sav_criterion}
Let $\mathcal{E}$ be a symmetric Dirichlet form on $L^2(X,\mm)$. Assume that there is a nondecreasing sequence $(X_n)_{n\ge1}$ of open subsets of $X$ such that:
\begin{enumerate}[label=(\alph*),ref=(\alph*),leftmargin=*]
    \item\label{item:crit_nest} $\mathrm{Cap}(X_n)<\infty$ and $\overline{X}_n\subseteq X_{n+1}$ for every $n$, and $(\overline{X}_n)_{n\ge1}$ is an $\mathcal{E}$-nest;
    \item\label{item:crit_compact} $\langle\mathrm{QR}.1'\rangle$: for all $n,m\ge1$ there is a compact set $K_{n,m}\subseteq X$ with $\mathrm{Cap}(X_n\setminus K_{n,m})\le1/m$;
    \item\label{item:crit_cont} $\langle\mathrm{QR}.2'\rangle$: the space of elements of $\V$ admitting a continuous representative is dense in $\V$ and separates the points of $X$.
\end{enumerate}
Then $\mathcal{E}$ is quasi-regular.
\end{lemma}

\begin{proposition}\label{prop:qr_facts}
Let $\mathcal{E}$ be a quasi-regular symmetric Dirichlet form on $L^2(X,\mm)$. Then:
\begin{enumerate}[label=(\roman*),ref=(\roman*),leftmargin=*]
    \item\label{item:qr_rep} every $u\in\V$ admits an $\mathcal{E}$-quasi-continuous representative $\tilde u$, unique up to $\mathcal{E}$-polar sets;
    \item\label{item:qr_measure} for every linear continuous $\ell\colon\V\to\R$ with $\ell(u)\geq0$ for all $u\in\V_+$, there is a unique $\sigma$-finite positive Borel measure $\mu$ on $X$, not charging $\mathcal{E}$-polar sets, such that $\tilde u\in L^1(\mu)$ and $\ell(u)=\int_X\tilde u\diff\mu$ for every $u\in\V$.
\end{enumerate}
\end{proposition}

\begin{proof}
\ref{item:qr_rep} is \cite[Rem.~1.3.9~(ii)]{CF12}, which holds for every quasi-regular Dirichlet form on a $\sigma$-finite measure space.

\ref{item:qr_measure}. The bilinear form $\mathcal{E}_1(u,v)\coloneqq\mathcal{E}(u,v)+\int_Xuv\diff\mm$ is the scalar product of $\V$, so by the Riesz representation theorem there is $u_\ell\in\V$ with $\mathcal{E}_1(u_\ell,v)=\ell(v)$ for every $v\in\V$; in particular, $\mathcal{E}_1(u_\ell,v)\geq0$ for every $v\in\V_+$. We claim that $u_\ell$ is $1$-excessive, that is, $u_\ell\geq0$ and $e^{-t}\PP_tu_\ell\leq u_\ell$ $\mm$-a.e.\ for every $t>0$.

For $u\in\V$, the Markov property gives $u^\pm\in\V$ and $\mathcal{E}(\abs{u})\leq\mathcal{E}(u)$; expanding $\abs{u}=u^++u^-$ and $u=u^+-u^-$, this reads $\mathcal{E}(u^+,u^-)\leq0$. Since $u_\ell^+u_\ell^-=0$, we get
\[
    0\leq\mathcal{E}_1(u_\ell,u_\ell^-)=\mathcal{E}(u_\ell^+,u_\ell^-)-\mathcal{E}_1(u_\ell^-)\leq-\mathcal{E}_1(u_\ell^-),
\]
so $u_\ell^-=0$. Fix now $v\in L^2(\mm)$ with $v\geq0$, and set $\phi(s)\coloneqq e^{-s}\int_Xu_\ell\,\PP_sv\diff\mm$. By Remark~\ref{rem:semigroup}, $\phi$ is continuous on $[0,\infty)$ and of class $C^1$ on $(0,\infty)$, and for $s>0$ we have $\PP_sv\in D(\mathsf{L})$ and $\PP_sv\geq0$, so that
\[
    \phi'(s)=e^{-s}\int_Xu_\ell\,(\mathsf{L}\PP_sv-\PP_sv)\diff\mm=-e^{-s}\,\mathcal{E}_1(u_\ell,\PP_sv)\leq0 .
\]
Hence, by the symmetry of $\PP_t$, $\int_Xe^{-t}\PP_tu_\ell\,v\diff\mm=\phi(t)\leq\phi(0)=\int_Xu_\ell\,v\diff\mm$, and the claim follows since $v$ is arbitrary.

Since $\mathcal{E}$ is quasi-regular and $u_\ell\in\V$ is $1$-excessive, the transfer method of \cite[Ch.~VI, Prop.~2.1]{MR92} gives a unique $\sigma$-finite positive Borel measure $\mu$ on $X$, not charging $\mathcal{E}$-polar sets, such that $\mathcal{E}_1(u_\ell,v)=\int_X\tilde v\diff\mu$ for every $v\in\V_+$. For $v\in\V$, the functions $(\tilde v)^\pm$ are quasi-continuous representatives of $v^\pm\in\V_+$, hence coincide q.e., and so $\mu$-a.e., with $\widetilde{v^\pm}$ by the uniqueness in \ref{item:qr_rep}. Therefore $\int_X(\tilde v)^\pm\diff\mu=\ell(v^\pm)<\infty$, so $\tilde v\in L^1(\mu)$ and $\int_X\tilde v\diff\mu=\ell(v^+)-\ell(v^-)=\ell(v)$.

This is the argument of \cite[(2.7) and Prop.~2.5]{Sav14}, where the statement is given under the standing assumptions of \cite[Section~2]{Sav14}: strong locality, existence of a carr\'e du champ and mass preservation of the semigroup. The proof above uses none of them.
\end{proof}

\begin{remark}\label{rem:qc}
Let $\mathcal{E}$ be quasi-regular, let $u\in\V$, and let $(F_k)_k$ be an $\mathcal{E}$-nest such that $\tilde u|_{F_k}$ is continuous for every $k$. Continuous functions are quasi-continuous with respect to the trivial nest $F_k=X$. By the uniqueness in Proposition~\ref{prop:qr_facts}\ref{item:qr_rep}:
\begin{enumerate}[label=(\roman*),ref=(\roman*),leftmargin=*]
    \item\label{item:qc_cont} if $u$ is continuous, then $\tilde u=u$ q.e.; hence $\ell(u)=\int_Xu\diff\mu$ for $\ell,\mu$ as in Proposition~\ref{prop:qr_facts}\ref{item:qr_measure};
    \item\label{item:qc_prod} if $\zeta\in C(X)\cap L^\infty(\mm)$ and $\zeta u\in\V$, then $(\zeta u)^\sim=\zeta\tilde u$ q.e.;
    \item\label{item:qc_comp} if $\eta\colon\R\to\R$ is Lipschitz with $\eta(0)=0$ and $\eta\circ u\in\V$, then $(\eta\circ u)^\sim=\eta\circ\tilde u$ q.e.
\end{enumerate}
Indeed, in \ref{item:qc_prod} and \ref{item:qc_comp} the right-hand sides are representatives of $\zeta u$ and $\eta\circ u$ whose restrictions to each $F_k$ are continuous.
\end{remark}

\begin{proposition}[The Cheeger energy as a Dirichlet form]\label{prop:identification}
Assume \textbf{(S)}, and let $\mathcal{E}(f)\coloneqq\int_X\abs{Df}^2\diff\mm$ for $f\in W^{1,2}(X)$ and $\mathcal{E}(f)\coloneqq+\infty$ for $f\in L^2(\mm)\setminus W^{1,2}(X)$. Then $\mathcal{E}=2\Ch$ is a strongly local symmetric Dirichlet form on $L^2(X,\mm)$ with domain $\V=W^{1,2}(X)$, it admits the carr\'e du champ $\Gamma(f,g)=\nabla f\cdot\nabla g$, its generator is $\Delta$, and its Markov semigroup is $(h_t)_{t\geq0}$.
\end{proposition}

\begin{proof}
By \textbf{(S)} and Lemma~\ref{lem:locally_finite}, $(X,\dd)$ is Polish and $\mm$ is $\sigma$-finite with full support.

\emph{Dirichlet form.} $\mathcal{E}$ is a quadratic form by \ref{item:IH}. It is lower semicontinuous on $L^2(\mm)$: if $f_n\to f$ in $L^2(\mm)$ and $\liminf_n\mathcal{E}(f_n)<\infty$, we extract a subsequence realising the $\liminf$, converging $\mm$-a.e., and along which $\abs{Df_n}\rightharpoonup G$ weakly in $L^2(\mm)$; by Proposition~\ref{prop:GP_calculus}\ref{item:calc_rules}, $G$ is a weak upper gradient of $f$, hence $f\in W^{1,2}(X)$ and $\mathcal{E}(f)\le\int_XG^2\diff\mm\le\liminf_n\mathcal{E}(f_n)$. For the Markov property, there is nothing to prove if $\mathcal{E}(f)=\infty$; if $f\in W^{1,2}(X)$ and $\eta$ is $1$-Lipschitz with $\eta(0)=0$, then $\abs{\eta\circ f}\le\abs{f}$ and $\abs{D(\eta\circ f)}\le\abs{Df}$ $\mm$-a.e.\ by the chain rule. The domain $W^{1,2}(X)$ is dense in $L^2(\mm)$ by Lemma~\ref{lem:LIPfm}.

\emph{Carr\'e du champ and strong locality.} By \ref{item:IH} and Proposition~\ref{prop:GP_calculus}\ref{item:calc_rules}, $(f,g) \mapsto \nabla f \cdot \nabla g$ is continuous, symmetric and bilinear, and satisfies the Leibniz rule; expanding with the Leibniz rule gives $2\int_X \Gamma(f,g)\varphi \diff\mm = \mathcal{E}(f,g\varphi) + \mathcal{E}(g,f\varphi) - \mathcal{E}(fg,\varphi)$ for $f,g,\varphi \in \V_\infty$. Let $f,g \in W^{1,2}(X)$ with $(f+a)g = 0$ $\mm$-a.e.\ for some $a \in \R$. Then $f = -a$ $\mm$-a.e.\ on $\{g \neq 0\}$, so $\abs{Df} = 0$ $\mm$-a.e.\ on $\{g \neq 0\}$ by locality, while $\abs{Dg} = 0$ $\mm$-a.e.\ on $\{g = 0\}$. Hence $\abs{\nabla f \cdot \nabla g} \leq \abs{Df}\abs{Dg} = 0$ $\mm$-a.e.\ and $\mathcal{E}(f,g) = 0$.

\emph{Generator and semigroup.} By Proposition~\ref{prop:GP_calculus}\ref{item:calc_lap}, $\int_Xg\,\Delta f\diff\mm=-\mathcal{E}(f,g)$ for $f\in D(\Delta)$ and $g\in W^{1,2}(X)$, so $\Delta$ is the generator of $\mathcal{E}$. Fix $f\in L^2(\mm)$ and set $u(t)\coloneqq h_tf$, $v(t)\coloneqq\PP_tf$. Both curves belong to $C([0,\infty);L^2(\mm))\cap C^1((0,\infty);L^2(\mm))$, take values in $D(\Delta)$ for $t>0$, start from $f$, and satisfy $u'=\Delta u$, $v'=\Delta v$ on $(0,\infty)$, by Proposition~\ref{prop:GP_calculus}\ref{item:calc_heat} and Remark~\ref{rem:semigroup}. Hence $z\coloneqq u-v$ satisfies
\[
    \frac{\diff}{\diff t}\norm{z(t)}_{L^2(\mm)}^2=2\int_Xz(t)\,\Delta z(t)\diff\mm=-2\int_X\abs{Dz(t)}^2\diff\mm\le0\qquad\text{for } t>0,
\]
and since $z$ is continuous on $[0,\infty)$ with $z(0)=0$, we get $z\equiv0$.
\end{proof}

\section{Proof of Theorem~\ref{thm:main}}\label{sec:proof}

Throughout this section we assume \textbf{(S)} and use the notation of Section~\ref{sec:dirichlet} with the identifications of Proposition~\ref{prop:identification}: $\V=W^{1,2}(X)$, $\mathcal{E}(u,v)=\int_X\G(u,v)\diff\mm$ with $\G(u,v)=\nabla u\cdot\nabla v$, $\G(u)=\abs{Du}^2$, $\mathsf{L}=\Delta$ and $\PP_t=h_t$. The proof is organised as follows:
\begin{gather*}
    L^2\text{-}\BE
    \implies \text{Quasi-regularity (Lemma~\ref{lem:QR})}\\
    \implies \text{Improved Bochner inequality (Proposition~\ref{prop:improved})}
    \implies L^1\text{-}\BE \;\text{ (Proposition~\ref{prop:L1BE})}\\
    \implies \text{Power-Harnack (Lemma~\ref{lem:power_harnack})}
    \implies \text{Theorem~\ref{thm:main}}.
\end{gather*}

\subsection{Quasi-regularity}\label{subsec:QR}

We now show that, under \textbf{(S)}, the Cheeger energy is a quasi-regular Dirichlet form. Quasi-regularity is what makes Proposition~\ref{prop:qr_facts} available: in Subsection~\ref{subsec:gamma2} it is used to represent the nonnegative functional $\GG^K[u]$ by a measure, and to evaluate it on Lipschitz functions. We verify Savar\'e's criterion (Lemma~\ref{lem:Sav_criterion}) on the open sets $X_n$ where finitely many elements of a dense sequence of test functions are not too small. Since $(X,\dd)$ need not be locally compact, the main point is condition~\ref{item:crit_compact}: compact sets of almost full capacity in $X_n$ are obtained by removing from $\overline{X}_n$ an open set of small capacity, outside of which $\overline{X}_n$ is totally bounded.

\begin{lemma}[Quasi-regularity]\label{lem:QR}
Assume \textbf{(S)}. Then $\mathcal{E}=2\Ch$ is a quasi-regular Dirichlet form on $L^2(X,\mm)$.
\end{lemma}

\begin{proof}
We verify the conditions of Lemma~\ref{lem:Sav_criterion}.

\emph{Step 1: condition~\ref{item:crit_cont}.} By Lemma~\ref{lem:GP}\ref{item:GP_dense}, $\Test\subseteq\V\cap\LIP(X)$ is dense in $\V$. Given $x\neq y$, let $B_r(x)$ be admissible with $2r<\dd(x,y)$; the function $\beta_{x,r}$ of Lemma~\ref{lem:LIPfm}\ref{item:LIPfm_bump} belongs to $\V\cap C(X)$ and $\beta_{x,r}(x)=1\neq0=\beta_{x,r}(y)$.

\emph{Step 2: condition~\ref{item:crit_nest}.} By Lemma~\ref{lem:sep_cap}, $\V$ is separable, so there are test functions $f_j\in\Test$, $j\geq1$, forming a dense subset of $\V$. Set
\[
    X_n\coloneqq\bigcup_{j=1}^n\Bigl\{\abs{f_j}>\frac{1}{n+1}\Bigr\},\qquad n\geq1 .
\]
Each $X_n$ is open, the sequence is nondecreasing, and by continuity of the $f_j$
\[
    \overline{X}_n\subseteq\bigcup_{j=1}^n\Bigl\{\abs{f_j}\geq\frac1{n+1}\Bigr\}\subseteq\bigcup_{j=1}^{n+1}\Bigl\{\abs{f_j}>\frac1{n+2}\Bigr\}=X_{n+1}.
\]
Since $(n+1)\abs{f_j}\in\V$ is at least $1$ on $\{\abs{f_j}>(n+1)^{-1}\}$, Lemma~\ref{lem:sep_cap} gives
\[
    \mathrm{Cap}(X_n)\leq(n+1)^2\sum_{j\leq n}\norm{f_j}_\V^2<\infty .
\]
To prove that $(\overline{X}_n)_n$ is a nest, fix $j$ and, for $n\geq j$, let $\eta_n(s)\coloneqq\mathrm{sgn}(s)\bigl(\abs{s}-\frac1{n+1}\bigr)^+$. Then $\eta_n\circ f_j\in\V$ vanishes outside $X_n$, so $\eta_n\circ f_j\in\V_{\overline{X}_n}$; moreover $\abs{f_j-\eta_n\circ f_j}\leq\abs{f_j}\wedge\frac1{n+1}$ and, by the chain rule and locality, $\abs{D(f_j-\eta_n\circ f_j)}=\abs{Df_j}\chi_{\{0<\abs{f_j}\leq(n+1)^{-1}\}}$ $\mm$-a.e. By dominated convergence $\eta_n\circ f_j\to f_j$ in $\V$, so $\bigcup_n\V_{\overline{X}_n}$ is dense in $\V$.

\emph{Step 3: condition~\ref{item:crit_compact}.} Fix $n\geq1$ and set
\[
    \Psi_n\coloneqq1\wedge\Bigl(2(n+1)\max_{1\leq j\leq n}\Bigl(\abs{f_j}-\frac{1}{2(n+1)}\Bigr)^+\Bigr),\qquad\Omega_n\coloneqq\bigcup_{j=1}^n\Bigl\{\abs{f_j}\geq\frac1{2(n+1)}\Bigr\}.
\]
$\Psi_n$ is Lipschitz, $0\leq\Psi_n\leq1$, $\Psi_n=1$ on $\overline{X}_n$, $\Psi_n=0$ outside $\Omega_n$, and, by Chebyshev's inequality,
\[
    \mm(\Omega_n)\leq4(n+1)^2\sum_{j\leq n}\norm{f_j}_{L^2(\mm)}^2<\infty .
\] Let $(x_i)_{i\geq1}$ be dense in $X$ and, for $k\geq1$,
\[
    w_k\coloneqq\Psi_n\wedge\min_{1\leq i\leq k}\dd(\cdot,x_i).
\]
Each $w_k$ is Lipschitz with $\mathrm{Lip}(w_k)\leq L_n\coloneqq\max\{\mathrm{Lip}(\Psi_n),1\}$, vanishes outside the closed set $\Omega_n$ and satisfies $0\leq w_k\leq\Psi_n\leq1$, so $w_k\in\LIPfm(X)$ and, by Lemma~\ref{lem:LIPfm}\ref{item:LIPfm_sob}, $\norm{w_k}_\V^2\leq(1+L_n^2)\mm(\Omega_n)$. The sequence $(w_k)_k$ is nonincreasing and $w_k\to0$ pointwise, since $(x_i)$ is dense; hence $w_k\to0$ in $L^2(\mm)$ and $w_k\rightharpoonup0$ weakly in $\V$. By the Banach--Saks theorem there is an increasing sequence $(k_i)_i$ such that $v_h\coloneqq h^{-1}\sum_{i=1}^hw_{k_i}\to0$ strongly in $\V$; we extract $(v_{h_l})_{l\geq1}$ with $\norm{v_{h_l}}_\V\leq2^{-2l}$. The functions $v_{h_l}$ are continuous, so $G_m\coloneqq\bigcup_{l\geq m}\{v_{h_l}>2^{-l}\}$ is open; the series $W_m\coloneqq\sum_{l\geq m}2^lv_{h_l}$ converges in $\V$ with $\norm{W_m}_\V\leq2^{-m+1}$, and $W_m\geq1$ $\mm$-a.e.\ on $G_m$. Hence the closed set $F_m\coloneqq X\setminus G_m$ satisfies $\mathrm{Cap}(X\setminus F_m)\leq2^{-2m+2}$. On $F_m$ we have $w_{k_{h_l}}\leq v_{h_l}\leq2^{-l}$ for $l\geq m$, because $(w_k)_k$ is nonincreasing. Given $\delta\in(0,1)$, choose $l\geq m$ with $2^{-l}<\delta$ and set $p\coloneqq k_{h_l}$: since $\Psi_n=1$ on $\overline{X}_n$, every $x\in\overline{X}_n\cap F_m$ satisfies $\min_{i\leq p}\dd(x,x_i)=w_p(x)<\delta$, that is, $\overline{X}_n\cap F_m\subseteq\bigcup_{i\leq p}B_\delta(x_i)$. Thus $\tilde K_{n,m}\coloneqq\overline{X}_n\cap F_m$ is closed and totally bounded, hence compact by completeness, and $\mathrm{Cap}(X_n\setminus\tilde K_{n,m})\leq\mathrm{Cap}(X\setminus F_m)\leq2^{-2m+2}$. The sets $K_{n,m}\coloneqq\tilde K_{n,m+2}$ satisfy $\mathrm{Cap}(X_n\setminus K_{n,m})\leq2^{-2m-2}<1/m$.
\end{proof}

\subsection{Test functions and the heat flow}\label{subsec:test}

In this subsection we collect the properties of test functions used in Subsections~\ref{subsec:gamma2}--\ref{subsec:L1BE}. Lemma~\ref{lem:approx} shows that every $u\in\V_\infty$ is approximated in $\V$ by test functions, with uniform $L^\infty$ bounds and preserving the sign; this allows us to define $\Gamma_2$ on test functions only and to test the resulting inequalities against general functions in $\V_\infty$. Lemma~\ref{lem:test_invariance} shows that $\Test$ is invariant under the heat flow, with a Lipschitz bound coming from \eqref{eq:BE_quadratic} and \ref{item:SL}. Lemma~\ref{lem:time_regularity} gives the time regularity of $r\mapsto\G(h_rf)$ as an $L^1(\mm)$-valued map, needed to differentiate the function $A$ of Subsection~\ref{subsec:L1BE}.

\begin{lemma}[Approximation by test functions]\label{lem:approx}
Assume \textbf{(S)}. Let $(\varphi_k)_{k \geq 1} \subseteq C_c^\infty(0,\infty)$ be nonnegative with $\int_0^\infty \varphi_k(s)\diff s = 1$ and $\mathrm{supp}(\varphi_k) \subseteq (0,1/k)$, and let $u \in \V_\infty$. Then $h_{\varphi_k}u \in \Test$, $\abs{h_{\varphi_k}u} \leq \norm{u}_{L^\infty(\mm)}$ $\mm$-a.e., $h_{\varphi_k}u \geq 0$ $\mm$-a.e.\ whenever $u \geq 0$ $\mm$-a.e., and $h_{\varphi_k}u \to u$ in $\V$ as $k \to \infty$. Consequently, $\Gamma(h_{\varphi_k}u) \to \Gamma(u)$ and $\Gamma(h_{\varphi_k}u,v) \to \Gamma(u,v)$ in $L^1(\mm)$ for every $v \in \V$.
\end{lemma}

\begin{proof}
The fact that $h_{\varphi_k}u$ belongs to $\Test$ is Lemma~\ref{lem:GP}\ref{item:GP_test}, and the two pointwise bounds follow from the weak maximum principle applied to $\pm u$ and $\norm{u}_{L^\infty(\mm)} \mp u$. For $s \downarrow 0$ we have $h_s u \to u$ in $L^2(\mm)$ and $\Ch(h_s u) \leq \Ch(u)$, so the lower semicontinuity of $\Ch$ gives $\Ch(h_s u) \to \Ch(u)$; since $\V$ is a Hilbert space, $h_s u \to u$ in $\V$. As $s \mapsto h_s u$ is continuous from $[0,\infty)$ to $\V$,
\[
    \norm{h_{\varphi_k}u - u}_{\V} = \Bigl\lVert \int_0^\infty (h_s u - u)\varphi_k(s)\diff s \Bigr\rVert_{\V} \leq \sup_{s \in (0,1/k)} \norm{h_s u - u}_{\V} \longrightarrow 0 .
\]
The $L^1$ convergences follow from the $\mm$-a.e.\ inequalities
\[
    \abs{\Gamma(w) - \Gamma(z)} \leq \Gamma(w - z)^{1/2}\bigl(\Gamma(w)^{1/2} + \Gamma(z)^{1/2}\bigr),\qquad \abs{\Gamma(w,v) - \Gamma(z,v)} \leq \Gamma(w - z)^{1/2}\Gamma(v)^{1/2},
\]
integrated with the Cauchy--Schwarz inequality.
\end{proof}

\begin{lemma}\label{lem:test_invariance}
Assume \textbf{(S)}. Let $f\in\Test$ and $r\geq0$. Then $h_rf\in\Test$, $\Delta h_rf=h_r\Delta f$, and $\G(h_rf)\leq e^{-2Kr}\mathrm{Lip}(f)^2$ $\mm$-a.e.
\end{lemma}

\begin{proof}
By \eqref{eq:BE_quadratic}, the bound $\abs{Df}\leq\mathrm{Lip}(f)$ (which follows from Lemma~\ref{lem:LIPfm}\ref{item:LIPfm_sob} by locality) and the weak maximum principle,
\[
    \G(h_rf)\leq e^{-2Kr}h_r\bigl(\G(f)\bigr)\leq e^{-2Kr}\mathrm{Lip}(f)^2 ,
\]
so $h_rf$ has a Lipschitz representative by \ref{item:SL}. The weak maximum principle gives $\norm{h_rf}_{L^\infty(\mm)}\leq\norm{f}_{L^\infty(\mm)}$. By Proposition~\ref{prop:GP_calculus}\ref{item:calc_heat}, $h_rf\in D(\Delta)$ with $\Delta h_rf=h_r\Delta f$, which belongs to $\V$ because $\Ch(h_rv)\leq\Ch(v)$, and to $L^\infty(\mm)$ by the weak maximum principle.
\end{proof}

\begin{lemma}\label{lem:time_regularity}
Assume \textbf{(S)}. Let $f\in\Test$, set $f_r\coloneqq h_rf$ and $G(r)\coloneqq\G(f_r)$ for $r\geq0$, and let $\psi\in C^2([0,\infty))$ with $\psi(0)=0$ and $\psi',\psi''$ bounded. Then:
\begin{enumerate}[label=(\roman*),ref=(\roman*),leftmargin=*]
    \item\label{item:tr_f} $r\mapsto f_r$ belongs to $C^1([0,\infty);\V)$, with derivative $\Delta f_r=h_r\Delta f$;
    \item\label{item:tr_G} $G\in C^1([0,\infty);L^1(\mm))$ with $G'(r)=2\G(f_r,\Delta f_r)$, and $\sup_{r\in[0,T]}\norm{G(r)}_{L^\infty(\mm)}<\infty$ for every $T>0$;
    \item\label{item:tr_psiG} $\psi\circ G\in C^1([0,\infty);L^1(\mm))\cap C([0,\infty);L^2(\mm))$ with $(\psi\circ G)'(r)=\psi'(G(r))\,G'(r)$.
\end{enumerate}
\end{lemma}

\begin{proof}
\ref{item:tr_f}. Set $v\coloneqq\Delta f\in\V$. For $\sigma,\tau\geq0$ we have $\norm{h_{\sigma+\tau}v-h_\sigma v}_\V\leq\norm{h_\tau v-v}_\V$, by the $L^2$-contractivity of $h_\sigma$ and $\Ch(h_\sigma w)\leq\Ch(w)$, and $h_\tau v\to v$ in $\V$ as $\tau\downarrow0$, as shown in the proof of Lemma~\ref{lem:approx}. Hence $\sigma\mapsto h_\sigma v$ is uniformly continuous from $[0,\infty)$ to $\V$. Since $f_r-f_{r'}=\int_{r'}^rh_\sigma v\diff\sigma$ in $L^2(\mm)$ for $0\leq r'<r$ and the integrand is continuous with values in $\V$, the identity holds in $\V$, and \ref{item:tr_f} follows.

\ref{item:tr_G}. This follows from \ref{item:tr_f} and the continuity of $\G\colon\V\times\V\to L^1(\mm)$; the $L^\infty$ bound is Lemma~\ref{lem:test_invariance}.

\ref{item:tr_psiG}. Fix $T>0$ and $M\coloneqq\sup_{[0,T]}\norm{G}_{L^\infty(\mm)}$. Since $\norm{w}_{L^2(\mm)}^2\leq\norm{w}_{L^\infty(\mm)}\norm{w}_{L^1(\mm)}$, $G$ is continuous with values in $L^2(\mm)$ on $[0,T]$, and so is $\psi\circ G$, $\psi$ being Lipschitz with $\psi(0)=0$. For $r,r+\tau\in[0,T]$, write $\Delta_\tau G\coloneqq G(r+\tau)-G(r)$ and $\Theta_\tau\coloneqq\int_0^1\psi'(G(r)+\theta\Delta_\tau G)\diff\theta$, so that
\[
    \frac{\psi(G(r+\tau))-\psi(G(r))}{\tau}-\psi'(G(r))G'(r)
    =\Theta_\tau\Bigl(\frac{\Delta_\tau G}{\tau}-G'(r)\Bigr)+\bigl(\Theta_\tau-\psi'(G(r))\bigr)G'(r).
\]
The first term tends to $0$ in $L^1(\mm)$ because $\abs{\Theta_\tau}\leq\norm{\psi'}_\infty$. For the second,
\[
    \abs[\big]{\Theta_\tau-\psi'(G(r))}\leq\min\bigl\{2\norm{\psi'}_\infty,\norm{\psi''}_\infty\abs{\Delta_\tau G}\bigr\}
\]
tends to $0$ in $\mm$-measure, and dominated convergence along subsequences gives convergence to $0$ in $L^1(\mm)$. The continuity of $r\mapsto\psi'(G(r))G'(r)$ in $L^1(\mm)$ follows in the same way.
\end{proof}

\subsection{The \texorpdfstring{$\Gamma_2$}{Gamma2} functional}\label{subsec:gamma2}

\begin{definition}\label{def:gamma2}
For $u,v\in\Test$ and $\varphi\in\V$ we set
\[
    \GG[u,v](\varphi)\coloneqq-\frac12\,\mathcal{E}\bigl(\G(u,v),\varphi\bigr)-\frac12\int_X\bigl(\G(u,\Delta v)+\G(v,\Delta u)\bigr)\varphi\diff\mm,
\]
$\GG^K[u,v](\varphi)\coloneqq\GG[u,v](\varphi)-K\int_X\G(u,v)\varphi\diff\mm$, $\GG[u]\coloneqq\GG[u,u]$ and $\GG^K[u]\coloneqq\GG^K[u,u]$. For $u,v,w\in\Test$ we set
\[
    H[u](v,w)\coloneqq\frac12\Bigl(\G\bigl(\G(u,v),w\bigr)+\G\bigl(\G(u,w),v\bigr)-\G\bigl(u,\G(v,w)\bigr)\Bigr).
\]
\end{definition}

\begin{remark}\label{rem:gamma2_wd}
By Lemma~\ref{lem:GP}\ref{item:GP_sq} and polarisation, $\G(u,v)\in\V_\infty$ for $u,v\in\Test$, so $\mathcal{E}(\G(u,v),\cdot)$ is continuous on $\V$; moreover $\G(u,\Delta v)\in L^2(\mm)$, since $\abs{Du}\in L^\infty(\mm)$ and $\Delta v\in\V$. Hence $\GG[u,v]$ and $\GG^K[u,v]$ are continuous linear functionals on $\V$, symmetric and bilinear in $(u,v)$, and $\GG^K[u,v]=\frac14\bigl(\GG^K[u+v]-\GG^K[u-v]\bigr)$. Similarly, $H[u](v,w)\in L^2(\mm)$ is symmetric in $(v,w)$, and
\begin{equation}\label{eq:H_gradient}
    \G\bigl(\G(u),w\bigr)=2H[u](u,w)\qquad\mm\text{-a.e.}
\end{equation}
\end{remark}

\begin{lemma}\label{lem:gamma2_pos}
Assume \textbf{(S)}. Then $\GG^K[u](\varphi)\geq0$ for every $u\in\Test$ and every $\varphi\in\V_+$.
\end{lemma}

\begin{proof}
For $\varphi\in\Testp$, Proposition~\ref{prop:GP_calculus}\ref{item:calc_lap} gives $\frac12\int_X\G(u)\Delta\varphi\diff\mm=-\frac12\mathcal{E}(\G(u),\varphi)$, so \ref{item:BE} reads $\GG^K[u](\varphi)\geq0$. The general case follows from Lemma~\ref{lem:GP}\ref{item:GP_dense} and the continuity of $\GG^K[u]$ on $\V$.
\end{proof}

\begin{lemma}\label{lem:product}
Assume \textbf{(S)}. Let $u,v,w\in\Test$ and $\varphi\in\V_\infty$. Then
\begin{equation}\label{eq:product_gamma2}
    \GG[u,vw](\varphi)=\GG[u,w](v\varphi)+\GG[u,v](w\varphi)+2\int_XH[u](v,w)\,\varphi\diff\mm
\end{equation}
and
\begin{equation}\label{eq:product_H}
    H[w^2](w,w)=2w\,H[w](w,w)+2\,\G(w)^2\qquad\mm\text{-a.e.}
\end{equation}
\end{lemma}

\begin{proof}
All products below belong to $\V_\infty$ by the Leibniz rule and Lemma~\ref{lem:GP}\ref{item:GP_sq}--\ref{item:GP_alg}. For $a\in\Test$ and $b\in\V_\infty$, the Leibniz rule and Proposition~\ref{prop:GP_calculus}\ref{item:calc_lap} give
\[
    \mathcal{E}(a\,b,\varphi)=\mathcal{E}(b,a\varphi)-2\int_X\varphi\,\G(a,b)\diff\mm-\int_Xb\,\varphi\,\Delta a\diff\mm .
\]
Since $\G(u,vw)=v\,\G(u,w)+w\,\G(u,v)$, applying this with $(a,b)=(v,\G(u,w))$ and $(a,b)=(w,\G(u,v))$ yields
\[
\begin{split}
    \mathcal{E}\bigl(\G(u,vw),\varphi\bigr)={}&\mathcal{E}\bigl(\G(u,w),v\varphi\bigr)+\mathcal{E}\bigl(\G(u,v),w\varphi\bigr)-2\int_X\varphi\bigl(\G(v,\G(u,w))+\G(w,\G(u,v))\bigr)\diff\mm\\
    &-\int_X\varphi\bigl(\G(u,w)\Delta v+\G(u,v)\Delta w\bigr)\diff\mm .
\end{split}
\]
By Proposition~\ref{prop:GP_calculus}\ref{item:calc_lap} and the Leibniz rule,
\[
    \G\bigl(u,\Delta(vw)\bigr)=v\,\G(u,\Delta w)+w\,\G(u,\Delta v)+\G(u,v)\Delta w+\G(u,w)\Delta v+2\,\G\bigl(u,\G(v,w)\bigr),
\]
and $\G(vw,\Delta u)=v\,\G(w,\Delta u)+w\,\G(v,\Delta u)$. Inserting these identities into the definition of $\GG[u,vw](\varphi)$, the terms containing $\Delta v$ and $\Delta w$ without a gradient cancel, and the remaining terms give \eqref{eq:product_gamma2}. For \eqref{eq:product_H}, the Leibniz rule gives
\begin{gather*}
    \G(w^2,w)=2w\,\G(w),\qquad\G\bigl(w^2,\G(w)\bigr)=2w\,\G\bigl(w,\G(w)\bigr),\\
    \G\bigl(2w\,\G(w),w\bigr)=2\,\G(w)^2+2w\,\G\bigl(\G(w),w\bigr),
\end{gather*}
so $H[w^2](w,w)=\G(\G(w^2,w),w)-\frac12\G(w^2,\G(w))=2\,\G(w)^2+w\,\G(\G(w),w)$, and \eqref{eq:H_gradient} concludes.
\end{proof}

In the smooth setting, $\GG^K[u](\varphi)=\int_X(\Gamma_2(u)-K\G(u))\varphi\diff\mm$ for a function $\Gamma_2(u)$; here $\GG^K[u]$ is only a functional on $\V$, and the next lemma shows that, being nonnegative by \ref{item:BE}, it is given by a measure. This is where quasi-regularity enters the proof.

\begin{lemma}\label{lem:repr}
Assume \textbf{(S)}. For every $u\in\Test$ there is a unique positive Borel measure $\mu_u$ on $X$, not charging $\mathcal{E}$-polar sets, such that $\GG^K[u](\varphi)=\int_X\tilde\varphi\diff\mu_u$ for every $\varphi\in\V$. Moreover, $\mu_u(B)<\infty$ for every admissible ball $B$, and $\GG^K[u](\varphi)=\int_X\varphi\diff\mu_u$ for every continuous $\varphi\in\V$.
\end{lemma}

\begin{proof}
Set $\ell\coloneqq\GG^K[u]$. By Remark~\ref{rem:gamma2_wd}, $\ell$ is a continuous linear functional on $\V$, and by Lemma~\ref{lem:gamma2_pos}, $\ell(\varphi)\geq0$ for every $\varphi\in\V_+$. By Lemma~\ref{lem:QR}, $\mathcal{E}$ is quasi-regular, so Proposition~\ref{prop:qr_facts}\ref{item:qr_measure} gives a unique positive Borel measure $\mu_u$, not charging $\mathcal{E}$-polar sets, with $\ell(\varphi)=\int_X\tilde\varphi\diff\mu_u$ for every $\varphi\in\V$. If $\varphi\in\V$ is continuous, then $\tilde\varphi=\varphi$ q.e.\ by Remark~\ref{rem:qc}\ref{item:qc_cont}, hence $\mu_u$-a.e., and $\ell(\varphi)=\int_X\varphi\diff\mu_u$. If $B=B_r(x)$ is admissible, the function $\beta\coloneqq\bigl(1-\dd(\cdot,B)/r\bigr)^+$ is continuous, belongs to $\V$ by Lemma~\ref{lem:LIPfm} and satisfies $\beta\geq\chi_B$, so $\mu_u(B)\leq\int_X\beta\diff\mu_u=\GG^K[u](\beta)<\infty$.
\end{proof}

We write $\mu_u=\rho_u\mm+\mu_u^s$ for the Lebesgue decomposition of $\mu_u$ with respect to $\mm$, with $\rho_u\geq0$ and $\mu_u^s\perp\mm$.

\begin{lemma}\label{lem:test_measures}
Assume \textbf{(S)}. Let $\Lambda_1,\Lambda_2$ be positive Borel measures on $X$, finite on admissible balls, such that $\int_X\varphi\diff\Lambda_1\geq\int_X\varphi\diff\Lambda_2$ for every nonnegative $\varphi\in\LIPb(X)$ vanishing outside some admissible ball. Then $\Lambda_1\geq\Lambda_2$.
\end{lemma}

\begin{proof}
Fix an admissible ball $B$; then $\Lambda\coloneqq(\Lambda_1-\Lambda_2)\llcorner B$ is a finite signed Borel measure. For an open set $U\subseteq B$, the functions $\varphi_k\coloneqq\min\{1,k\,\dd(\cdot,X\setminus U)\}$ are Lipschitz, vanish outside $U$ and satisfy $0\leq\varphi_k\uparrow\chi_U$, so $\Lambda(U)\geq0$ by dominated convergence. Let $B=P\cup N$ be a Hahn decomposition for $\Lambda$ and $\delta>0$. The finite Borel measure $\abs{\Lambda}$ on the metric space $X$ is outer regular, so there is an open set $V\supseteq N$ with $\abs{\Lambda}(V\setminus N)<\delta$; then $U\coloneqq V\cap B$ is open and $0\leq\Lambda(U)\leq\Lambda(N)+\delta$. Hence $\Lambda(N)=0$, i.e.\ $\Lambda_1\llcorner B\geq\Lambda_2\llcorner B$. Since $X$ is covered by countably many admissible balls (Lemma~\ref{lem:locally_finite}), $\Lambda_1\geq\Lambda_2$.
\end{proof}

\subsection{The \texorpdfstring{$L^1$}{L1} Bakry--\'Emery estimate}\label{subsec:L1BE}

We follow the self-improvement argument of Bakry \cite{Bak86}, in the form of \cite[Section~3]{Sav14}. This gives the improved Bochner inequality $\G(\G(f))\leq4\rho_f\G(f)$ (Proposition~\ref{prop:improved}). Then, to obtain the $L^1$ Bakry--\'Emery estimate, we differentiate
\[
    A(s)\coloneqq\int_Xh_s\varphi\,\psi\bigl(\G(h_{t-s}f)\bigr)\diff\mm,\qquad s\in[0,t],\qquad\psi(r)\coloneqq\sqrt{r+\varepsilon}-\sqrt\varepsilon,
\]
for nonnegative $\varphi\in L^1(\mm)\cap L^\infty(\mm)$: by Proposition~\ref{prop:improved}, $A'\geq KA$ up to an error vanishing as $\varepsilon\downarrow0$. Lemma~\ref{lem:product_derivative} is the elementary differentiation rule used in this step.

\begin{proposition}[Improved Bochner inequality]\label{prop:improved}
Assume \textbf{(S)} and let $f\in\Test$. Then
\begin{equation}\label{eq:improved}
    \G\bigl(\G(f)\bigr)\leq4\rho_f\,\G(f)\qquad\mm\text{-a.e.}
\end{equation}
\end{proposition}

\begin{proof}
\emph{Step 1: expansion.} Fix $w\in\Test$ and $\lambda,c\in\R$. By Lemma~\ref{lem:GP}\ref{item:GP_alg}, $u\coloneqq\lambda f+w^2-2cw$ belongs to $\Test$. We write $H_f\coloneqq H[f](w,w)$ and $H_w\coloneqq H[w](w,w)$. We claim that, for every $\varphi\in\V_\infty$,
\begin{equation}\label{eq:expansion}
\begin{split}
    \GG^K[u](\varphi)={}&\lambda^2\,\GG^K[f](\varphi)+4\lambda\,\GG^K[f,w]\bigl((w-c)\varphi\bigr)+4\,\GG^K[w]\bigl((w-c)^2\varphi\bigr)\\
    &+\int_X\Bigl(4\lambda H_f+8(w-c)H_w+4\,\G(w)^2\Bigr)\varphi\diff\mm .
\end{split}
\end{equation}
Indeed, with $g\coloneqq w^2-2cw$ we have $\GG[u]=\lambda^2\GG[f]+2\lambda\GG[f,g]+\GG[g]$. By \eqref{eq:product_gamma2},
\begin{gather*}
    \GG[f,g](\varphi)=2\,\GG[f,w]\bigl((w-c)\varphi\bigr)+2\int_XH_f\,\varphi\diff\mm,\\
    \GG[w,w^2](\varphi)=2\,\GG[w](w\varphi)+2\int_XH_w\varphi\diff\mm,
\end{gather*}
and, using also the symmetry of $\GG$ and \eqref{eq:product_H},
\[
\begin{split}
    \GG[w^2](\varphi)&=2\,\GG[w^2,w](w\varphi)+2\int_XH[w^2](w,w)\,\varphi\diff\mm\\
    &=4\,\GG[w](w^2\varphi)+\int_X\bigl(8wH_w+4\,\G(w)^2\bigr)\varphi\diff\mm .
\end{split}
\]
Hence
\[
\begin{split}
    \GG[g](\varphi)&=\GG[w^2](\varphi)-4c\,\GG[w^2,w](\varphi)+4c^2\,\GG[w](\varphi)\\
    &=4\,\GG[w]\bigl((w-c)^2\varphi\bigr)+\int_X\bigl(8(w-c)H_w+4\,\G(w)^2\bigr)\varphi\diff\mm .
\end{split}
\]
Finally, $\G(u)=\lambda^2\G(f)+4\lambda(w-c)\G(f,w)+4(w-c)^2\G(w)$ by the Leibniz rule, and subtracting $K\int_X\G(u)\varphi\diff\mm$ yields \eqref{eq:expansion}.

\emph{Step 2: an inequality between measures.} Set $\nu_{f,w}\coloneqq\frac14(\mu_{f+w}-\mu_{f-w})$, so that $\GG^K[f,w](\varphi)=\int_X\varphi\diff\nu_{f,w}$ for every continuous $\varphi\in\V$, and
\[
    \Lambda_{\lambda,c}\coloneqq\lambda^2\mu_f+4\lambda(w-c)\,\nu_{f,w}+4(w-c)^2\mu_w+\bigl(4\lambda H_f+8(w-c)H_w+4\,\G(w)^2\bigr)\mm .
\]
Splitting each term into its positive and negative parts, $\Lambda_{\lambda,c}$ is the difference of two positive Borel measures which are finite on admissible balls: $\mu_f$, $\mu_w$, $\mu_{f\pm w}$ are finite on admissible balls by Lemma~\ref{lem:repr}, $w$ is bounded, and $H_f$, $H_w$, $\G(w)^2$ belong to $L^1(\mm)$. Let $\varphi\in\LIPb(X)$ vanish outside an admissible ball $B_r(x)$; then $\mathrm{supp}(\varphi)\subseteq\overline{B}_r(x)$, so $\varphi\in\LIPfm(X)$, and $\varphi$, $(w-c)\varphi$, $(w-c)^2\varphi$ are continuous elements of $\V$. By Lemma~\ref{lem:repr}, the right-hand side of \eqref{eq:expansion} equals $\int_X\varphi\diff\Lambda_{\lambda,c}$, and by Lemmas~\ref{lem:gamma2_pos} and \ref{lem:test_measures}, $\Lambda_{\lambda,c}\geq0$ for every $\lambda,c\in\R$.

\emph{Step 3: densities.} Fix an admissible ball $B$, let $\sigma\coloneqq(\mm+\mu_f+\mu_w+\mu_{f+w}+\mu_{f-w})\llcorner B$, a finite measure, and write $M$, $F$, $W$, $N$ for the densities with respect to $\sigma$ of $\mm\llcorner B$, $\mu_f\llcorner B$, $\mu_w\llcorner B$, $\nu_{f,w}\llcorner B$. Fix Borel representatives of $H_f$, $H_w$, $\G(w)$; changing them on an $\mm$-null set $Z$ changes the products with $M$ only on the $\sigma$-null set $Z\cap\{M>0\}$. By Step~2, for each $(\lambda,c)\in\Q^2$ the density
\[
    D_{\lambda,c}\coloneqq\lambda^2F+4\lambda(w-c)N+4(w-c)^2W+\bigl(4\lambda H_f+8(w-c)H_w+4\,\G(w)^2\bigr)M
\]
is nonnegative outside a $\sigma$-null set $E_{\lambda,c}$. Outside $E\coloneqq\bigcup_{(\lambda,c)\in\Q^2}E_{\lambda,c}$, the polynomial $(\lambda,c)\mapsto D_{\lambda,c}(x)$ is nonnegative on $\Q^2$, hence on $\R^2$, and the choice $c=w(x)$ gives
\[
    \lambda^2F(x)+4\lambda H_f(x)M(x)+4\,\G(w)(x)^2M(x)\geq0\qquad\text{for all }\lambda\in\R\text{ and }x\in B\setminus E .
\]
Since $\mm(\{M=0\}\cap B)=0$, and since $\mu_f=\rho_f\mm+\mu_f^s$ with $\mu_f^s$ concentrated on an $\mm$-null set (see after Lemma~\ref{lem:repr}), so that $F=\rho_fM$ $\mm$-a.e.\ on $B$, dividing by $M$ we get $\lambda^2\rho_f+4\lambda H_f+4\,\G(w)^2\geq0$ for all $\lambda\in\R$, $\mm$-a.e.\ on $B$. Since $X$ is covered by countably many admissible balls,
\begin{equation}\label{eq:hessian_quadratic}
    \abs{H[f](w,w)}\leq\sqrt{\rho_f}\,\G(w)\qquad\mm\text{-a.e.\ in }X.
\end{equation}

\emph{Step 4: polarisation.} Let $v,w\in\Test$. Applying \eqref{eq:hessian_quadratic} to $sv\pm s^{-1}w$ for $s\in\Q\setminus\{0\}$, we get $\mm$-a.e.
\[
\begin{split}
    \abs{H[f](v,w)}&=\frac14\abs[\big]{H[f](sv+s^{-1}w,sv+s^{-1}w)-H[f](sv-s^{-1}w,sv-s^{-1}w)}\\
    &\leq\frac14\sqrt{\rho_f}\bigl(\G(sv+s^{-1}w)+\G(sv-s^{-1}w)\bigr)=\frac12\sqrt{\rho_f}\bigl(s^2\G(v)+s^{-2}\G(w)\bigr),
\end{split}
\]
and the infimum over $s\in\Q\setminus\{0\}$ gives $\abs{H[f](v,w)}\leq\sqrt{\rho_f}\abs{Dv}\abs{Dw}$ $\mm$-a.e. With $v=f$ and \eqref{eq:H_gradient},
\begin{equation}\label{eq:pre_improved}
    \abs[\big]{\G\bigl(\G(f),w\bigr)}\leq2\sqrt{\rho_f}\,\abs{Df}\,\abs{Dw}\qquad\mm\text{-a.e., for every }w\in\Test .
\end{equation}

\emph{Step 5: approximation.} By Lemma~\ref{lem:GP}\ref{item:GP_sq} and Lemma~\ref{lem:approx}, there are $w_k\in\Test$ with $w_k\to\G(f)$ in $\V$. Then $\G(\G(f),w_k)\to\G(\G(f))$ in $L^1(\mm)$ and $\abs{Dw_k}\to\abs{D\G(f)}$ in $L^2(\mm)$, hence $\mm$-a.e.\ along a subsequence, and \eqref{eq:pre_improved} gives $\G(\G(f))\leq2\sqrt{\rho_f}\abs{Df}\sqrt{\G(\G(f))}$ $\mm$-a.e., which is \eqref{eq:improved}.
\end{proof}

\begin{lemma}\label{lem:product_derivative}
Let $T>0$ and let $a,b\colon[0,T]\to L^2(\mm)$ satisfy
\[
    a\in C([0,T];L^2(\mm))\cap C^1((0,T];L^2(\mm)),\qquad b\in C([0,T];L^2(\mm))\cap C^1([0,T];L^1(\mm)),
\]
and $\sup_{[0,T]}\norm{a}_{L^\infty(\mm)}<\infty$. Then $A(s)\coloneqq\int_Xa(s)b(s)\diff\mm$ belongs to $C([0,T])\cap C^1((0,T))$, and
\[
    A'(s)=\int_X\bigl(a'(s)b(s)+a(s)b'(s)\bigr)\diff\mm\qquad\text{for } s\in(0,T).
\]
\end{lemma}

\begin{proof}
Continuity follows from the continuity of the pairing on $L^2(\mm)\times L^2(\mm)$. For $s\in(0,T)$ and small $\tau\neq0$,
\[
    \frac{A(s+\tau)-A(s)}{\tau}=\int_X\frac{a(s+\tau)-a(s)}{\tau}\,b(s+\tau)\diff\mm+\int_Xa(s)\,\frac{b(s+\tau)-b(s)}{\tau}\diff\mm ;
\]
the first integral converges by convergence in $L^2(\mm)$ of both factors, the second because $a(s)\in L^\infty(\mm)$ and the difference quotient converges in $L^1(\mm)$. The continuity of $A'$ follows in the same way.
\end{proof}

\begin{proposition}[$L^1$ Bakry--\'Emery estimate]\label{prop:L1BE}
Assume \textbf{(S)}. Then for every $f \in \V$ and every $t \geq 0$,
\begin{equation}\label{eq:BE_1}
    \abs{D h_t f} \leq e^{-Kt} h_t\bigl(\abs{Df}\bigr) \quad \mm\text{-a.e.\ in } X.
\end{equation}
\end{proposition}

\begin{proof}
\emph{Step 1: setup.} Let first $f\in\Test$ and $t>0$, fix $\varepsilon>0$ and
\[
    \psi(r)\coloneqq\sqrt{r+\varepsilon}-\sqrt\varepsilon,\qquad\psi'(r)=\frac{1}{2\sqrt{r+\varepsilon}},\qquad\psi''(r)=-\frac{1}{4(r+\varepsilon)^{3/2}} ,
\]
and fix $\varphi\in L^1(\mm)\cap L^\infty(\mm)$ with $\varphi\geq0$. Set $a(s)\coloneqq h_s\varphi$ for $s\in[0,t]$. By Proposition~\ref{prop:GP_calculus}\ref{item:calc_heat}, $a\in C([0,t];L^2(\mm))\cap C^1((0,t];L^2(\mm))$ with $a'(s)=\Delta a(s)$, $0\leq a(s)\leq\norm{\varphi}_{L^\infty(\mm)}$ and $\int_Xa(s)\diff\mm\leq\norm{\varphi}_{L^1(\mm)}$. By Lemma~\ref{lem:time_regularity}\ref{item:tr_psiG}, the map $s\mapsto\psi(\G(h_{t-s}f))$ belongs to $C([0,t];L^2(\mm))\cap C^1([0,t];L^1(\mm))$, with derivative
\[
    -2\psi'\bigl(\G(h_{t-s}f)\bigr)\,\G(h_{t-s}f,\Delta h_{t-s}f).
\]
Let
\[
    A(s)\coloneqq\int_Xa(s)\,\psi\bigl(\G(h_{t-s}f)\bigr)\diff\mm,\qquad s\in[0,t].
\]
By Lemma~\ref{lem:product_derivative}, applied to $a$ and $s\mapsto\psi(\G(h_{t-s}f))$, $A\in C([0,t])\cap C^1((0,t))$ and, writing $g\coloneqq h_{t-s}f\in\Test$ (Lemma~\ref{lem:test_invariance}) and $a\coloneqq a(s)$,
\[
    A'(s)=\int_X\Delta a\;\psi\bigl(\G(g)\bigr)\diff\mm-2\int_Xa\,\psi'\bigl(\G(g)\bigr)\,\G(g,\Delta g)\diff\mm\qquad\text{for } s\in(0,t).
\]

\emph{Step 2: the derivative as a $\Gamma_2$ term.} Let $\Psi\in C^1(\R)$ be the extension of $\psi'$ given by $\Psi(r)\coloneqq\psi'(0)+\psi''(0)r$ for $r<0$; $\Psi$ is Lipschitz and positive, and $\Psi\geq\psi'(\norm{\G(g)}_{L^\infty(\mm)})$ on $(-\infty,\norm{\G(g)}_{L^\infty(\mm)}]$. By the chain rule, $\psi(\G(g))\in\V$ and
\[
    \chi\coloneqq\Psi\bigl(\G(g)\bigr)\,a=\Psi(0)\,a+\bigl(\Psi(\G(g))-\Psi(0)\bigr)a\in\V_\infty ,\qquad\nabla\chi=\psi'\bigl(\G(g)\bigr)\nabla a+a\,\psi''\bigl(\G(g)\bigr)\nabla\G(g).
\]
Hence, by Proposition~\ref{prop:GP_calculus}\ref{item:calc_lap},
\[
    \int_X\Delta a\;\psi\bigl(\G(g)\bigr)\diff\mm=-\int_X\psi'\bigl(\G(g)\bigr)\nabla a\cdot\nabla\G(g)\diff\mm=-\mathcal{E}\bigl(\G(g),\chi\bigr)+\int_Xa\,\psi''\bigl(\G(g)\bigr)\G\bigl(\G(g)\bigr)\diff\mm,
\]
and therefore
\begin{equation}\label{eq:A_prime}
    A'(s)=2\,\GG[g](\chi)+\int_Xa\,\psi''\bigl(\G(g)\bigr)\G\bigl(\G(g)\bigr)\diff\mm .
\end{equation}

\emph{Step 3: lower bound.} By Lemma~\ref{lem:GP}\ref{item:GP_Lip}, $a$ has a Lipschitz representative $\hat a$, which is nonnegative everywhere because $\mathrm{supp}(\mm)=X$. By Remark~\ref{rem:qc}\ref{item:qc_prod}--\ref{item:qc_comp} applied to the decomposition of $\chi$ above, $\tilde\chi=\Psi\bigl(\G(g)^\sim\bigr)\,\hat a\geq0$ q.e., hence $\mu_g$-a.e. By Lemma~\ref{lem:repr} and $\mu_g\geq\rho_g\mm$,
\[
    \GG[g](\chi)=\int_X\tilde\chi\diff\mu_g+K\int_X\G(g)\chi\diff\mm\geq\int_X\chi\,\rho_g\diff\mm+K\int_X\G(g)\chi\diff\mm ;
\]
in particular $a\rho_g\in L^1(\mm)$, since $\Psi(\G(g))\geq\psi'(\norm{\G(g)}_{L^\infty(\mm)})>0$. 
Since $\G(g)\geq0$ and $\Psi=\psi'$ on $[0,\infty)$, we have $\chi=\psi'(\G(g))\,a$, and inserting the last inequality into \eqref{eq:A_prime} gives
\[
    A'(s)\geq\int_Xa\Bigl(2\rho_g\,\psi'\bigl(\G(g)\bigr)+2K\,\G(g)\,\psi'\bigl(\G(g)\bigr)+\psi''\bigl(\G(g)\bigr)\,\G\bigl(\G(g)\bigr)\Bigr)\diff\mm .
\]
By Proposition~\ref{prop:improved} applied to $g\in\Test$, $\G(\G(g))\leq4\rho_g\G(g)$ $\mm$-a.e., and since $\psi''<0$,
\[
    \psi''\bigl(\G(g)\bigr)\,\G\bigl(\G(g)\bigr)\geq4\rho_g\,\G(g)\,\psi''\bigl(\G(g)\bigr)\qquad\mm\text{-a.e.}
\]
Hence
\[
    A'(s)\geq\int_Xa\Bigl(2\rho_g\Bigl(\psi'\bigl(\G(g)\bigr)+2\,\G(g)\,\psi''\bigl(\G(g)\bigr)\Bigr)+2K\,\G(g)\,\psi'\bigl(\G(g)\bigr)\Bigr)\diff\mm .
\]
For $r\geq0$ we have
\[
    \psi'(r)+2r\psi''(r)=\frac{\varepsilon}{2(r+\varepsilon)^{3/2}}\geq0,\qquad2r\psi'(r)=\frac{r}{\sqrt{r+\varepsilon}} .
\]
Since $a\geq0$ and $\rho_g\geq0$, the term containing $\rho_g$ is nonnegative, and
\[
    A'(s)\geq K\int_Xa\,\frac{\G(g)}{\sqrt{\G(g)+\varepsilon}}\diff\mm .
\]
Since, for $r\geq0$,
\[
    0\leq\frac{r}{\sqrt{r+\varepsilon}}-\psi(r)=\sqrt\varepsilon-\frac{\varepsilon}{\sqrt{r+\varepsilon}}\leq\sqrt\varepsilon ,
\]
and $\int_Xa\diff\mm\leq\norm{\varphi}_{L^1(\mm)}$, the right-hand side differs from $KA(s)=K\int_Xa\,\psi(\G(g))\diff\mm$ by at most $\abs{K}\sqrt\varepsilon\,\norm{\varphi}_{L^1(\mm)}$. Hence
\[
    A'(s)\geq K A(s)-\abs{K}\sqrt\varepsilon\,\norm{\varphi}_{L^1(\mm)}\qquad\text{for every }s\in(0,t).
\]

\emph{Step 4: integration.} Since $A$ is continuous on $[0,t]$ and of class $C^1$ on $(0,t)$, Gronwall's lemma gives
\[
    A(0)\leq e^{-Kt}A(t)+C_{K,t}\sqrt\varepsilon\,\norm{\varphi}_{L^1(\mm)},\qquad C_{K,t}\coloneqq\abs{K}\int_0^te^{-Ks}\diff s ,
\]
where $A(0)=\int_X\varphi\,\psi(\G(h_tf))\diff\mm$ and $A(t)=\int_Xh_t\varphi\,\psi(\G(f))\diff\mm$. Letting $\varepsilon\downarrow0$ by monotone convergence, since $\psi(r)=r\bigl(\sqrt{r+\varepsilon}+\sqrt\varepsilon\bigr)^{-1}$ increases to $\sqrt r$ as $\varepsilon\downarrow0$, and using the self-adjointness of $h_t$,
\[
    \int_X\varphi\,\abs{Dh_tf}\diff\mm\leq e^{-Kt}\int_Xh_t\varphi\,\abs{Df}\diff\mm=e^{-Kt}\int_X\varphi\,h_t\bigl(\abs{Df}\bigr)\diff\mm .
\]
Choosing $\varphi$ as the indicator of $\{\abs{Dh_tf}>e^{-Kt}h_t(\abs{Df})\}\cap E_n$, where $(E_n)_{n\in\N}$ are Borel sets of finite measure with $E_n\uparrow X$ (Lemma~\ref{lem:locally_finite}), and letting $n\to\infty$ gives \eqref{eq:BE_1} for $f\in\Test$.

\emph{Step 5: density.} For $f\in\V$, take $f_k\in\Test$ with $f_k\to f$ in $\V$ (Lemma~\ref{lem:GP}\ref{item:GP_dense}). Then $\abs{Dh_tf_k}\to\abs{Dh_tf}$ and $h_t(\abs{Df_k})\to h_t(\abs{Df})$ in $L^2(\mm)$, and \eqref{eq:BE_1} passes to the limit $\mm$-a.e.\ along a subsequence. The case $t=0$ is trivial.
\end{proof}

\subsection{The Harnack inequality}\label{subsec:harnack}

\begin{lemma}[Power-type Harnack inequality]\label{lem:power_harnack}
Assume \textbf{(S)}. Let $f \in \LIPfm(X)$ with $0 \leq f \leq 1/2$ and let $t > 0$. Then $h_t f$ has a Lipschitz representative, still denoted by $h_t f$, and
\begin{equation}\label{eq:harnack_bound}
    \bigl(h_t f(x)\bigr)^2 \leq h_t f(y) \exp\Bigl( \frac{\dd(x,y)^2}{2 I_K(t)} \Bigr) \qquad \text{for all } x, y \in X,
\end{equation}
where $I_K(t)$ is given by \eqref{eq:cK}.
\end{lemma}

\begin{proof}
For $r \in (0,t]$ set $f_r \coloneqq h_r f$. Since $f\in\LIPfm(X)$, we have $f\in L^1(\mm)\cap L^\infty(\mm)\cap\V$ by Lemma~\ref{lem:LIPfm}. By positivity and the weak maximum principle, $0 \leq f_r \leq 1/2$. By Lemma~\ref{lem:GP}\ref{item:GP_Lip}, $f_r$ admits a Lipschitz representative; moreover $f_r \in D(\Delta) \cap L^1(\mm) \cap L^\infty(\mm)$. We always use these representatives; since $\mathrm{supp}(\mm) = X$, two continuous functions that agree $\mm$-a.e.\ agree everywhere.

Fix $\delta \in (0,1/4)$ and consider the smooth convex function $\phi_\delta \colon [0,1/2] \to [0,\infty)$,
\[
    \phi_\delta(r) \coloneqq (r+\delta)\log(r+\delta) - r\log\delta - \delta\log\delta .
\]
Since $\phi_\delta(0) = 0$ and $\phi_\delta''(r) = (r+\delta)^{-1} \in [1,\delta^{-1}]$ on $[0,1/2]$, by Proposition~\ref{prop:GP_calculus}\ref{item:calc_lap}, applied to a $C^2$ extension of $\phi_\delta$ to $\R$ with bounded second derivative, we have $\phi_\delta(f_{t-s}) \in D(\Delta) \cap L^2(\mm)$ for every $s \in [0,t)$, with
\begin{equation}\label{eq:lap_phi_delta}
    \Delta\bigl(\phi_\delta(f_{t-s})\bigr) = \phi_\delta'(f_{t-s})\Delta f_{t-s} + \frac{\abs{Df_{t-s}}^2}{f_{t-s} + \delta}.
\end{equation}
Fix $g \in L^2(\mm) \cap L^\infty(\mm)$ with $g \geq 0$, and set $\Phi_\delta(s) \coloneqq \int_X h_s g\,\phi_\delta(f_{t-s})\diff\mm$ for $s \in [0,t]$. The curve $s \mapsto h_s g$ belongs to $C([0,t];L^2(\mm)) \cap C^1((0,t];L^2(\mm))$ with derivative $\Delta h_s g$ and $h_s g \in D(\Delta)$ for $s > 0$. Since $\phi_\delta$ is of class $C^2$ on $[0,1/2]$ with $\phi_\delta(0) = 0$, the curve $s \mapsto \phi_\delta(f_{t-s})$ belongs to $C([0,t];L^2(\mm)) \cap C^1([0,t);L^2(\mm))$ with derivative $-\phi_\delta'(f_{t-s})\Delta f_{t-s}$. Hence $\Phi_\delta \in C([0,t]) \cap C^1((0,t))$, and differentiating on $(0,t)$, using \eqref{eq:lap_phi_delta} and the symmetry of $\Delta$ on $D(\Delta)$, we get
\begin{equation}\label{eq:Phi_delta_deriv}
    \Phi_\delta'(s) = \int_X h_s g\,\frac{\abs{Df_{t-s}}^2}{f_{t-s}+\delta}\diff\mm
    = \int_X g\,h_s\Bigl(\frac{\abs{Df_{t-s}}^2}{f_{t-s}+\delta}\Bigr)\diff\mm .
\end{equation}

For $\lambda\in\R$, expanding $\bigl((f_{t-s}+\delta)^{-1/2}\abs{Df_{t-s}}+\lambda(f_{t-s}+\delta)^{1/2}\bigr)^2\geq0$ and applying the positivity-preserving operator $h_s$, with $h_sf_{t-s}=f_t$ and $h_s\delta\leq\delta$, we obtain
\[
    h_s\Bigl(\frac{\abs{Df_{t-s}}^2}{f_{t-s}+\delta}\Bigr)+2\lambda\,h_s\bigl(\abs{Df_{t-s}}\bigr)+\lambda^2(f_t+\delta)\geq0\quad\mm\text{-a.e.}
\]
For each $\lambda$ this holds outside an $\mm$-null set $N_\lambda$. The set $N\coloneqq\bigcup_{\lambda\in\Q}N_\lambda$ is $\mm$-null, and by continuity in $\lambda$ the inequality holds at every point of $X\setminus N$ for every $\lambda\in\R$. At such points, the choice $\lambda=-h_s(\abs{Df_{t-s}})/(f_t+\delta)$ gives
\begin{equation}\label{eq:CS_hs}
    \bigl(h_s(\abs{Df_{t-s}})\bigr)^2\leq(f_t+\delta)\,h_s\Bigl(\frac{\abs{Df_{t-s}}^2}{f_{t-s}+\delta}\Bigr)\quad\mm\text{-a.e.}
\end{equation}
Then, applying Proposition~\ref{prop:L1BE} at time $s$ to $f_{t-s}\in\V$ and using \eqref{eq:CS_hs}, we get
\begin{equation}\label{eq:pointwise_grad}
    \abs{Df_t}^2\leq e^{-2Ks}\bigl(h_s(\abs{Df_{t-s}})\bigr)^2\leq e^{-2Ks}(f_t+\delta)\,h_s\Bigl(\frac{\abs{Df_{t-s}}^2}{f_{t-s}+\delta}\Bigr)\quad\mm\text{-a.e.}
\end{equation}
Dividing \eqref{eq:pointwise_grad} by $e^{-2Ks}(f_t+\delta)$ and inserting it into \eqref{eq:Phi_delta_deriv}, we obtain
\[
    \Phi_\delta'(s)\geq e^{2Ks}\int_Xg\,\frac{\abs{Df_t}^2}{f_t+\delta}\diff\mm\qquad\text{for every } s\in(0,t).
\]
Since $\Phi_\delta\in C([0,t])$, integrating over $s\in(0,t)$ and writing $\Phi_\delta(t)=\int_Xg\,h_t(\phi_\delta(f))\diff\mm$ by the self-adjointness of $h_t$, we get
\[
    I_K(t)\int_Xg\,\frac{\abs{Df_t}^2}{f_t+\delta}\diff\mm\leq\Phi_\delta(t)-\Phi_\delta(0)=\int_Xg\,\bigl(h_t(\phi_\delta(f))-\phi_\delta(f_t)\bigr)\diff\mm .
\]
Both $\abs{Df_t}^2/(f_t+\delta)\leq\delta^{-1}\abs{Df_t}^2$ and $h_t(\phi_\delta(f))-\phi_\delta(f_t)$ belong to $L^1(\mm)$, the latter because $\phi_\delta$ is Lipschitz on $[0,1/2]$ with $\phi_\delta(0)=0$ and $f,f_t\in L^1(\mm)$. Since $g\in L^2(\mm)\cap L^\infty(\mm)$ with $g\geq0$ is arbitrary, we obtain
\begin{equation}\label{eq:entropy_prod_bound}
    I_K(t)\,\frac{\abs{Df_t}^2}{f_t+\delta}\leq h_t\bigl(\phi_\delta(f)\bigr)-\phi_\delta(f_t)\quad\mm\text{-a.e.}
\end{equation}

We now bound the right-hand side of \eqref{eq:entropy_prod_bound} from above. The linear part of $\phi_\delta$ gives no contribution: writing
\[
    \phi_\delta(r)=\theta_\delta(r)-r\log\delta,\qquad\theta_\delta(r)\coloneqq(r+\delta)\log(r+\delta)-\delta\log\delta,
\]
the linearity of $h_t$ and $h_tf=f_t$ give
\[
    h_t\bigl(\phi_\delta(f)\bigr)-\phi_\delta(f_t)=h_t\bigl(\theta_\delta(f)\bigr)-\theta_\delta(f_t).
\]
Moreover $0<f+\delta\leq3/4<1$, so $(f+\delta)\log(f+\delta)\leq0$ and hence $\theta_\delta(f)\leq\theta_\delta(f)^+\leq-\delta\log\delta$ $\mm$-a.e. Since $\theta_\delta$ is Lipschitz on $[0,1/2]$ with $\theta_\delta(0)=0$, we have $\theta_\delta(f)^+\in L^2(\mm)$, and positivity together with the weak maximum principle gives
\[
    h_t\bigl(\theta_\delta(f)\bigr)\leq h_t\bigl(\theta_\delta(f)^+\bigr)\leq-\delta\log\delta .
\]
Therefore
\[
    h_t\bigl(\phi_\delta(f)\bigr)-\phi_\delta(f_t)\leq-\delta\log\delta-\theta_\delta(f_t)=(f_t+\delta)\log\frac{1}{f_t+\delta},
\]
and \eqref{eq:entropy_prod_bound} becomes
\begin{equation}\label{eq:fisher_bound}
    I_K(t)\,\abs{Df_t}^2\leq(f_t+\delta)^2\log\frac{1}{f_t+\delta}\quad\mm\text{-a.e.}
\end{equation}

Consider $\eta_\delta(r)\coloneqq\sqrt{-\log\delta}-\sqrt{-\log(r+\delta)}$ for $r\in[0,1/2]$. Since $r+\delta\in[\delta,3/4]$, we have $-\log(r+\delta)\geq\log(4/3)>0$, so $\eta_\delta$ is smooth on $[0,1/2]$ with $\eta_\delta(0)=0$ and
\[
    0<\eta_\delta'(r)=\frac{1}{2(r+\delta)\sqrt{-\log(r+\delta)}}\leq\frac{1}{2\delta\sqrt{\log(4/3)}} .
\]
Extending $\eta_\delta$ to a Lipschitz function on $\R$, the chain rule gives $w_\delta\coloneqq\eta_\delta(f_t)\in\V$ with $\abs{Dw_\delta}=\eta_\delta'(f_t)\abs{Df_t}$ $\mm$-a.e., and by \eqref{eq:fisher_bound}
\[
    \abs{Dw_\delta}^2=\frac{\abs{Df_t}^2}{4(f_t+\delta)^2\bigl(-\log(f_t+\delta)\bigr)}\leq\frac{1}{4I_K(t)}\quad\mm\text{-a.e.}
\]
By \ref{item:SL}, $w_\delta$ has a $(2\sqrt{I_K(t)})^{-1}$-Lipschitz representative, which coincides with the continuous function $\eta_\delta(f_t)$. Hence, for all $x,y\in X$,
\[
    \sqrt{-\log(f_t(y)+\delta)}\leq\sqrt{-\log(f_t(x)+\delta)}+\frac{\dd(x,y)}{2\sqrt{I_K(t)}} .
\]
Squaring with $(a+b)^2\leq2a^2+2b^2$ and exponentiating gives
\[
    (f_t(x)+\delta)^2\leq(f_t(y)+\delta)\exp\Bigl(\frac{\dd(x,y)^2}{2I_K(t)}\Bigr),
\]
and letting $\delta\downarrow0$ yields \eqref{eq:harnack_bound}.
\end{proof}

\subsection{Proof of Theorem~\ref{thm:main}}\label{subsec:proof_main}

\begin{proof}[Proof of Theorem~\ref{thm:main}]
The function $f$ is $(2r_0)^{-1}$-Lipschitz, $0\leq f\leq1/2$ and $\mathrm{supp}(f)\subseteq\overline{B}_{r_0}(x_0)$, so $f\in\LIPfm(X)$ and Lemma~\ref{lem:power_harnack} applies.

\emph{Non-degeneracy.} The heat flow is injective on $L^2(\mm)$. Indeed, let $g \in L^2(\mm)$ and $s > 0$ with $h_s g = 0$. By the semigroup property and the self-adjointness of $h_{s/2}$,
\[
    \norm{h_{s/2} g}_{L^2(\mm)}^2 = \int_X g\,h_{s/2}h_{s/2}g\diff\mm = \int_X g\,h_s g\diff\mm = 0 ,
\]
so $h_{s/2}g = 0$; by induction $h_{s2^{-k}}g = 0$ for every $k \in \N$. Since $\tau \mapsto h_\tau g$ is continuous from $[0,\infty)$ to $L^2(\mm)$, letting $k \to \infty$ gives $g = 0$. Now $f \neq 0$ in $L^2(\mm)$ because $\mathrm{supp}(\mm) = X$, hence $h_t f \neq 0$. Since $h_t f \geq 0$ is continuous, there is $\bar{x} \in X$ with $h_t f(\bar{x}) > 0$. Applying \eqref{eq:harnack_bound} with $(x,y) = (\bar{x},x_0)$ gives $h_t f(x_0) \geq \bigl(h_t f(\bar{x})\bigr)^2 \exp\bigl(-\dd(\bar{x},x_0)^2/(2I_K(t))\bigr) > 0$.

\emph{Integration.} Applying \eqref{eq:harnack_bound} with $(x,y) = (x_0,y)$, for every $y \in B_r(x_0)$ we get
\[
    h_t f(y) \geq \bigl(h_t f(x_0)\bigr)^2 \exp\Bigl(-\frac{r^2}{2I_K(t)}\Bigr).
\]
Integrating over $B_r(x_0)$ and using $h_t f \geq 0$ and the $L^1(\mm)$-contractivity of $h_t$,
\[
    \mm\bigl(B_r(x_0)\bigr)\bigl(h_t f(x_0)\bigr)^2 \exp\Bigl(-\frac{r^2}{2I_K(t)}\Bigr) \leq \int_{B_r(x_0)} h_t f \diff\mm \leq \norm{f}_{L^1(\mm)},
\]
which is \eqref{eq:main_growth}. Since every bounded set is contained in a ball centred at $x_0$, $\mm$ is finite on bounded sets; then $(X,\dd,\mm)$ satisfies the assumptions of Definition~\ref{def:RCD}, with \ref{item:G} given by \eqref{eq:main_growth} for $\bar x=x_0$, and it is an $\RCD(K,\infty)$ space.
\end{proof}

\subsection*{Acknowledgements}
The author thanks Nicola Gigli for carefully checking the arguments of this note. This work was supported by the Italian Ministry of University and Research (MUR) through the programme ``Fondo Italiano per la Scienza'' (FIS~3), project ``Modern challEnges in Geometric Analysis -- MEGA'' (PI: N.~Gigli), CUP G53C25000920001. Parts of the arguments in this note were developed, checked and typeset with the assistance of the large language models Claude Opus~5.5 (Anthropic) and Gemini~4 Argon (Google DeepMind). The author takes full responsibility for the content.
The author is a member of the GNAMPA group of INdAM.

\subsection*{Data availability}
Data sharing is not applicable to this article, as no datasets were generated or analysed during the current study.

\subsection*{Conflict of interest}
The author has no competing interests to declare that are relevant to the content of this article.


\end{document}